# Solutions with a uniform time of existence of a class of characteristic semi-linear wave equations near $\mathscr{S}^+$

Marcel Dossa
Département de Mathématiques
Faculté des Sciences, Université de Yaoundé I

Roger Tagne Wafo
Département de Mathématiques et Informatique
Faculté des Sciences, Université de Douala

July 24, 2013

## Abstract

We prove existence and uniqueness of solution of a class of semi-linear wave equations with initial data prescribed on the light-cone with vertex the origin of a Minkowski space-time. The nonlinear term is assumed to satisfy a nullity condition which guarantee that the neighborhood of the initial cone on which we obtain our solution does not shrink to zero as one approaches infinity. This result is applied to wave maps on Minkowski space-time $\mathbb{R}^{n+1}$ with $n \geq 3$.

# Contents

# 1 Introduction

Let $(\mathbb{R}^{n+1}_x, \eta_x)$ be the usual Minkowski space time with the global canonical coordinates system $(x^\mu)$. We denote by $\mathcal{C}^+_{a,x}$ the translated half cone of equation $x^0 = r + a$ where $a > 0$, $r^2 = \sum_{i=1}^n (x^i)^2$, $r \geq 0$ and by $\mathcal{Y}^+_{a,x}$ the interior of $\mathcal{C}^+_{a,x}$, that is the set of points $(x^\mu)$ such that $x^0 > r + a$. In this work, we are interested with the following characteristic semi-linear Cauchy problem

$$\begin{cases} \Box_{x,\eta_x} f = F(\cdot, f, \partial f) & \text{in } \mathcal{Y}^+_{a,x} \\ f = \varphi & \text{on } \mathcal{C}^+_{a,x} \end{cases} \qquad (1.1)$$

where $\eta_x = (\eta_{\alpha\beta})$ is the Minkowski metric on $\mathbb{R}^{n+1}_x$, $\eta = diag(-1, +1, \ldots, +1)$, $\Box_{x,\eta_x}$ the flat wave operator,

$$f = (f^I), \quad \partial f = \left(\frac{\partial f^I}{\partial x^\alpha}\right), \quad F = (F^I), \quad \alpha = 0, 1 \ldots, n, \ I = 1, \ldots, N \ ,$$

and

$$\varphi = (\varphi^I), \quad \text{the initial data prescribed on} \quad \mathcal{C}^+_{a,x} \ .$$

There exists in the literature a complete study (even in the quasi-linear case) of problem (1.1) near the tip of the cone $\mathcal{C}^+_{a,x}$, see the series of papers [3, 7, 9, 10] and the references therein; compare [19, 18, 15] for a very general treatment of Lipschitz initial data hypersurfaces for the linear wave equation. Under suitable conditions on the source term and/or on the initial data, in those papers, it is shown that, in the semi-linear or quasi-linear case, there exists a neighborhood of the tip of the initial cone in $\mathcal{Y}^+_{a,x}$ on which one can find a unique solution. As far as the global solution of (1.1) is concerned (i.e. existence of solutions in the entire interior of the initial cone), a lot remains to be done. It is well known that for an arbitrary nonlinear function $F$, in general it is not possible to solve globally or semi-globally this problem, that is, without restriction on $F$ and/or on the space dimension $n$, it is not possible to find a neighborhood of the whole half cone $\mathcal{C}^+_{a,x}$ on which one can get existence and uniqueness of solution of such problem. In [1], A. Cabet gave some examples of nonlinearities for which the solution develops singularities in finite time regardless the smallness and/or the smoothness of the initial data in the case $n = 1$. To the best of our knowledge, three types of nonlinearities have been considered so far, leading to global or semi-global solution of (1.1):

- In [13], M. Dossa and F. Touadera assume that the space dimension $n$ is odd and greater than or equal to 3, that the source term $F = F(f, \partial f)$ is such that $F(0,0) = F'(0,0) = 0$ and $F^{(2)}$ satisfies the null condition of S. Klainerman when $n = 3$. With these conditions on the nonlinear term and the space dimension $n$, it is shown that if the initial data prescribed on the light cone are small enough in some appropriate norms, then (1.1) has a global solution in the whole interior of the initial cone.

- In [11], the authors suppose that, the restriction to the initial cone of the functions $F(x^\mu, f(x^\mu), \partial f(x^\mu))$ is a linear function with respect to the restriction to same cone of the derivatives of the unknown function $f(x^\mu)$ with respect to $x^0$. With this hypothesis, they proved that there exists a neighborhood of the entire initial cone on which problem (1.1) has a unique solution. We notice that this result does not guarantee that the thickness of the obtained neighborhood does not shrink to zero as one approaches infinity.

- In [12, 16, 1], analogous characteristic Cauchy problem are considered with initial data specify on two intersecting smooth null hypersurfaces under some suitable null condition on $F$. The results of these last references combined with local existence results on a neighborhood of the tip of the cone $\mathcal{C}^+_{a,x}$ of [9, 8] can also permit to study problem (1.1) under the condition that $F$ is linear with respect to the derivatives of the unknown function in the normal direction of the initial cone $\mathcal{C}^+_{a,x}$. Indeed, assuming this, one succeed in concluding as in the previous case. We should point here that in the reference [1], it remains to fix a problem of regularity of initial data and of dependance of some constants used in the iterative scheme on $\lambda$. [1]

The difficulty here is due to the fact that, in the process of solving such problem, one needs to estimate the outgoing derivatives of the unknown function on the initial cone. The characteristic property of this cone does not allow to choose arbitrarily the first of these derivatives as it is the case in the classical Cauchy problem. In order to obtain global solution, we need to solve globally a nonlinear ordinary differential equation with a nonlinear part which is exactly $F$. In the third case we mentioned above, this equation is linear and thus can be globally solved on $\mathcal{C}^+_{a,x}$. We intend in this paper to show that there exists a future neighborhood not only of the entire null cone $\mathcal{C}^+_{a,x}$ but also by guaranteeing that the thickness of this neighborhood does not nullify when one reaches infinity, on which there exists a unique solution of (1.1). To do this, we shall impose on the function $F$ a hypothesis of nullity of the kind of [6], see hypothesis 4.21 of this reference. More precisely, we shall suppose that the function $F$ has a *uniform zero in* $(f, \partial f) = 0$ of order $r \geq 2$ which is related to the space dimension (regarless the fact that $n$ is odd or not) by the condition

$$n \geq 1 + \frac{4}{r-1} - 2\alpha , \tag{1.2}$$

and that the initial data $\varphi$ are in some weighted Sobolev spaces near the conformal infinity. The non negative real number $-\alpha$ is the exponent of the weight in our Sobolev's norms, which is chosen in order to control the singularities arising in the volume element from the gauge transformation we use. In the particular case $\alpha = -\frac{1}{2}$, the constraints (1.2) reads:

| $r =$ | 2 | 3, 4 | 5, 6, ... |
|---|---|---|---|
| $n \geq$ | 6 | 4 | 3 |

The strategy here will be based on the techniques of conformal method used in [4, 5] by P.T. Chruściel and R. T. Wafo in the case of the classical Cauchy problem, the method of iterative scheme introduced in [17] by A. Majda and repeated by A. Cabet in [1] and R. Racke in [20] and finally the method of local solution developed by M. Dossa in [7], see also [9] and [10].

## 2 Transformation of the system

### 2.1 The transformed wave equation

Let $\mathcal{C}_{0,x}$ be the light cone of Minkowski space time $\mathbb{R}^{n+1}_x$ of equation $(x^0)^2 = r^2$. We will denote by $\mathcal{C}^+_{0,x}$ and $\mathcal{C}^-_{0,x}$ the future and past light cone of the origin of coordinates respectively, by $\mathcal{Y}^+_{0,x}$ the interior of $\mathcal{C}^+_{0,x}$ and by $\mathcal{Y}^-_{0,x}$ the interior of $\mathcal{C}^-_{0,x}$. As in [6], we consider the map $\phi$ defined as:

$$\phi : \mathbb{R}^{n+1}_x \setminus \mathcal{C}_{0,x} \to \mathbb{R}^{n+1}_y \text{ by } x^\alpha \mapsto y^\alpha := \frac{x^\alpha}{\eta_{\lambda\mu} x^\lambda x^\mu} , \ \alpha = 0, 1, \ldots, n . \tag{2.1}$$

---

[1] In that reference, the definition of the surface element $dS' = e^{-\lambda \Psi_+} dS$ on the slices $N^-_u |_V \equiv [0, V] \times Y$ in the unnumbered equation after equation 4.2 page 2115, implies that the Sobolev constant $c'$ of equation 4.4 page 2116, depends on $\lambda$. In fact as it is said there, $c' = c_s e^{\lambda V}$ where $c_s$ is an universal Sobolev constant coming from the embedding $H^m(U) \hookrightarrow C^1(U)$, $U$ subset of $\mathbb{R}^n$ and $m > \frac{n}{2} + 1$. The consequence is that the constant $\tilde{c}_3(\rho)$ might depends exponentially on $\lambda$ and it will not be possible to choose $\lambda$ such that $\tilde{c}_3(\rho) - \lambda \tilde{c} \leq 0$ as stated there.

Note that $\phi(\mathbb{R}_x^{n+1} \setminus \mathcal{C}_{0,x})$ is a subset fo $\mathbb{R}_y^{n+1}$. Any of the sets defined in $\mathbb{R}_x^{n+1}$ has its counterpart in $\mathbb{R}_y^{n+1}$, we keep the same notations. The indices $x$ or $y$ will be used to indicate if the set under consideration is a subset of $\mathbb{R}_x^{n+1}$ or $\mathbb{R}_y^{n+1}$. As an example, the set $\mathcal{C}_{0,y}$ is the light cone with vertex the origin of coordinates in $\mathbb{R}_y^{n+1}$, its equation is given by $(y^0)^2 = \rho^2$ where $\rho^2 = \sum_{i=1}^{n}(y^i)^2$.

The conformal map $\phi$ $\left(ds_x^2 = (-\eta_{\alpha\beta} y^\alpha y^\beta)^{-2} ds_y^2\right)$ is a bijection from $\mathcal{Y}_{0,x}^+$ onto $\phi(\mathcal{Y}_{0,x}^+) = \mathcal{Y}_{0,y}^-$, with inverse

$$\phi^{-1} : y^\alpha \mapsto x^\alpha \quad \text{by} \quad x^\alpha := \frac{y^\alpha}{\eta_{\lambda\mu} y^\lambda y^\mu} \ . \tag{2.2}$$

$\phi$ is also a bijection from $\mathcal{Y}_{a,x}^+$ onto the relatively compact domain $\phi(\mathcal{Y}_{a,x}^+) = \mathcal{Y}_{-\frac{1}{a},y}^+ \cap \mathcal{Y}_{0,y}^-$ (see Figure 2 page 41) with the same inverse as in (2.2).

Further, by setting

$$\Omega = -\eta_{\alpha\beta} y^\alpha y^\beta \quad \text{and} \quad \hat{f} = \Omega^{-\frac{n-1}{2}} f \circ \phi^{-1} \ , \tag{2.3}$$

one obtains(the details of calculations can be found in [22]):

$$\frac{\partial f}{\partial x^\mu} \circ \phi^{-1}(y^\nu) = \Omega^{\frac{n-1}{2}} \left\{ (1-n) y_\mu \hat{f} - \Omega \frac{\partial \hat{f}}{\partial y^\mu} - 2 y_\mu y^\alpha \frac{\partial \hat{f}}{\partial y^\alpha} \right\} \ ; \tag{2.4}$$

and

$$\Box_{x,\eta} f = \Omega^{\frac{n+3}{2}} \Box_{y,\eta} \hat{f} \ . \tag{2.5}$$

Thus the right-hand side of equation (1.1) reads:

$$\begin{aligned} F\left(x^\nu, f(x^\nu), \partial_\mu f(x^\nu)\right) &= F\left(\phi^{-1}(y^\nu), f \circ \phi^{-1}(y^\nu), \partial_\mu f \circ \phi^{-1}(y^\nu)\right) \\ &= F\left(\phi^{-1}(y^\nu), \Omega^{\frac{n-1}{2}} \hat{f}, \Omega^{\frac{n-1}{2}} \left\{(1-n) y_\mu \hat{f} - \Omega \frac{\partial \hat{f}}{\partial y^\mu} - 2 y_\mu y^\alpha \frac{\partial \hat{f}}{\partial y^\alpha} \right\}\right) \\ &\equiv \widetilde{F}\left(y^\nu, \Omega^{\frac{n-1}{2}} \hat{f}, \Omega^{\frac{n-1}{2}} \frac{\partial \hat{f}}{\partial y^\mu}\right) \ . \end{aligned}$$

We obtain that under the coordinates transformation (2.1) and the rescaling (2.3), the wave equation (1.1) reads:

$$\begin{cases} \Box_{y,\eta} \hat{f} = \Omega^{-\frac{n+3}{2}} \widetilde{F}\left(y^\nu, \Omega^{\frac{n-1}{2}} \hat{f}, \Omega^{\frac{n-1}{2}} \frac{\partial \hat{f}}{\partial y^\mu}\right) & \text{in} \quad \phi(\mathcal{Y}_{a,x}^+) \\ \hat{f} = \hat{\varphi} \quad \text{on} \quad \phi(\mathcal{C}_{a,x}^+) \end{cases} \tag{2.6}$$

where

$$\hat{\varphi} = \left.\left(\Omega^{-\frac{n-1}{2}} f \circ \phi^{-1}\right)\right|_{\phi(\mathcal{C}_{a,x}^+)} \ .$$

**Remark** 2.1 To the system 2.6 we can apply the results of [10] to obtain that there exists a neighborhood which will be denoted by $V_{0,y}$ (see Figure 3 page 41) of the tip of the cone $\phi(\mathcal{C}_{a,x}^+)$ on which (2.6) has a unique smooth solution. We denote this local solution by $\hat{f}_0$.

## 2.2 Goursat problem associated to the transformed system

As in [2], we consider the Cauchy problem associated to the wave equation of system (2.6) with prescribed data on two truncated (such as to get rid of the tips) intersecting cones $\mathcal{C}^+ \subset \mathcal{C}_{-\frac{1}{a},y}^+ \cap \mathcal{Y}_{0,y}^-$ and $\mathcal{C}^- \subset \mathcal{C}_{\lambda,y}^- \cap \mathcal{Y}_{-\frac{1}{a},y}^+$, where $\lambda$ is a fixed parameter belonging to the interval $]-\frac{1}{a}, 0[$ sufficiently close to $-\frac{1}{a}$ such that $\mathcal{C}_{\lambda,y}^- \cap \mathcal{Y}_{-\frac{1}{a},y}^+$ intercepts $V_{0,y}$ (see Figures 3 and 4 pages 41 and 42):

$$\Box_{y,\eta}\hat{f} = \Omega^{-\frac{n+3}{2}}\widetilde{F}\left(y^\nu, \Omega^{\frac{n-1}{2}}\hat{f}, \Omega^{\frac{n-1}{2}}\frac{\partial \hat{f}}{\partial y^\mu}\right) \ ; \tag{2.7}$$

in the future neighborhood of $\mathcal{C}^+ \cup \mathcal{C}^-$ with initial data

$$\hat{f} = \hat{\varphi} \quad \text{on} \quad \mathcal{C}^+ \quad \text{and} \quad \hat{f} = \hat{f}_0 \quad \text{sur} \quad \mathcal{C}^- \ ; \tag{2.8}$$

where $\hat{f}_0$ is the smooth function given by Remark 2.1 in the neighborhood $V_{0,y}$ of the tip $(-\frac{1}{a}, 0)$. We will be concerned now in deriving a global process which solves (2.7)-(2.8). The next section is devoted to this goal.

## 3 Existence and uniqueness theorem

### 3.1 Second transformation

In the space $\mathbb{R}^{n+1}_y$ we consider now the spherical coordinates $(\tau, \rho, \theta)$ defined as:

$$\begin{cases} \tau = y^0, \\ \rho = \left(\sum_{i=1}^n (y^i)^2\right)^{1/2}, \\ y^i = \rho \omega^i(\theta), \ i = 1, \ldots, n \end{cases} \quad \text{with} \quad \begin{cases} \omega^1 = \cos\theta^1 \\ \omega^2 = \sin\theta^1 \cos\theta^2 \\ \omega^3 = \sin\theta^1 \sin\theta^2 \cos\theta^3 \\ \ldots \quad \ldots \\ \omega^{n-1} = \sin\theta^1 \sin\theta^2 \ldots \sin\theta^{n-2} \cos\theta^{n-1} \\ \omega^n = \sin\theta^1 \sin\theta^2 \ldots \sin\theta^{n-2} \sin\theta^{n-1} \end{cases}$$

where $0 < \theta^{n-1} < 2\pi$ and $0 < \theta^i < \pi$, $i = 1, 2, \ldots, n-2$. We set:

$$\begin{cases} x = \tau + \rho \leq 0 \\ y = \tau - \rho + \frac{1}{a} \geq 0 \end{cases} \quad \text{i.e.} \quad \begin{cases} \tau = \frac{1}{2}(y + x - \frac{1}{a}) \\ \rho = \frac{1}{2}(\frac{1}{a} + x - y) \end{cases} . \tag{3.1}$$

In the new coordinate system $(y, x, \theta)$, we have the identity

$$\Box_{y,\eta} = -4\partial_x \partial_y + \frac{n-1}{\rho}(\partial_x - \partial_y) + \frac{\Delta_{S^{n-1}}}{\rho^2} \tag{3.2}$$

where $\Delta_{S^{n-1}}$ is the Laplace-Beltrami operator on the sphere $S^{n-1}$ endowed with its canonical metric. From this identity, we deduce the new form of the transformed equation (2.6) with respect to the new coordinates system $z := (y, x, \theta)$:

$$\begin{cases} -4\partial_x \partial_y \hat{f} + \frac{n-1}{\rho}\left(\partial_x - \partial_y\right)\hat{f} + \frac{\Delta_{S^{n-1}}\hat{f}}{\rho^2} = \Omega^{-\frac{n+3}{2}}\widetilde{F}\left(z, \Omega^{\frac{n-1}{2}}\hat{f}, \Omega^{\frac{n-1}{2}}\frac{\partial \hat{f}}{\partial y^\mu}\right) \text{ in } \phi(\mathcal{V}^+_{a,x}) \\ \hat{f} = \hat{\varphi} \quad \text{on} \quad \phi(\mathcal{C}^+_{a,x}) \end{cases} \tag{3.3}$$

**Remark 3.1** We emphasise on the fact that $\Omega = -x(1/a - y)$ and $y^\mu \frac{\partial}{\partial y^\mu} = x\partial_x + (y - 1/a)\partial_y$. Thus by identity (2.4), we will suppose without further restriction on $F$ that when replacing the first order derivatives $\frac{\partial \hat{f}}{\partial y^\mu}$ in $\widetilde{F}$ by their value in terms of $\partial_y$, $\partial_x$, $\partial_A$ (we will write $\partial_A$ for $\partial_{\theta^A}$), any derivative $\partial_x \hat{f}$ comes with a pre-factor $x$.

**Remark 3.2** In the coordinates system $(\tau, \rho, \theta)$ the Minkowski metric reads:

$$\eta = -(d\tau)^2 + \sum_{i=1}^n (dy^i)^2 = -(d\tau)^2 + (d\rho)^2 + \rho^2 ds^2$$

where

$$\begin{aligned} ds^2 &= (d\theta^1)^2 + \sin^2\theta^1 (d\theta^2)^2 + \sin^2\theta^1 \sin^2\theta^2 (d\theta^3)^2 + \ldots \\ &\quad + \sin^2\theta^1 \sin^2\theta^2 \ldots \sin^2\theta^{n-2}(d\theta^{n-1})^2 \end{aligned} \tag{3.4}$$

thus
$$\eta = -(d\tau)^2 + (d\rho)^2 + \rho^2 h_{AB}^2 d\theta^A d\theta^B,$$
with $h_{AB} = 0$ if $A \neq B$ and $h_{AA}$, $A = 1, \ldots, n-1$ being defined by equation (3.4). The inverse metric is then given by
$$\eta^\sharp = -(\partial_\tau)^2 + (\partial_\rho)^2 + h^{AB}\partial_{\theta^A}\partial_{\theta^B} \quad \text{with} \quad h^{AB} = \begin{cases} 0 & \text{if } A \neq B \\ \frac{1}{\rho^2 h_{AA}} & \text{if } A = B \end{cases}.$$

**Remark** 3.3

- $\phi(\mathcal{Y}_{a,x}^+) = \mathcal{Y}_{-\frac{1}{a},y}^+ \cap \mathcal{Y}_{0,y}^- = \{(y,x,\theta) : -\frac{\sqrt{2}}{2a} \leq x \leq 0; \ 0 \leq y \leq \frac{\sqrt{2}}{2a}\};$

$$\phi(\mathcal{C}_{a,x}^+) = \mathcal{C}_{-\frac{1}{a},y}^+ \cap \mathcal{Y}_{0,y}^- = \{(y,x,\theta) : y = 0; \ -\frac{\sqrt{2}}{2a} \leq x \leq 0\}.$$

- $\frac{1}{a} + x - y = 0$ is equivalent to $\rho = 0$, thus the function $x \mapsto \frac{1}{\frac{1}{a} - x - y}$ is well defined as far as one does not reach $\{\rho = 0\}$ (which will be the case in the domain of interest).

- Setting $x_0 = \sqrt{2}/2\lambda$, $y_0 = \sqrt{2}/2(\lambda + 1/a)$ we have:
$$\mathcal{C}^- = \mathcal{C}_{\lambda,y}^- \cap \mathcal{Y}_{-\frac{1}{a},y}^+ = \{(y,x,\theta) : x = x_0, \ 0 \leq y \leq y_0\}.$$
and
$$\mathcal{C}^+ = \{(y,x,\theta) : y = 0, \ x_0 \leq x < 0\}.$$

## 3.2 Functional spaces

We intend in this section to describe the slices (see Figure 5 page 42) on which we will get our energy estimates. Let $z := (y,x,\theta)$, be a generic point and denote by $\mathcal{D}$ the set defined by $\mathcal{D} = [0, y_0] \times [x_0, 0[ \times \mathcal{O}$, where $\mathcal{O}$ is a subset of the unit sphere $\mathbb{S}^{n-1}$ of $\mathbb{R}^n$. For any $(u,v) \in [0, y_0] \times [x_0, 0[$, we set
$$\mathcal{D}_{u,v} = [0, u] \times [x_0, v] \times \mathcal{O},$$
$$\mathcal{C}_{u,v}^+ = \{u\} \times [x_0, v] \times \mathcal{O} = \{(y,x,\theta) : y = u; \ x_0 \leq x \leq v\}$$
and
$$\mathcal{C}_{u,v}^- = [0, u] \times \{v\} \times \mathcal{O} = \{(y,x,\theta) : 0 \leq y \leq u; \ x = v\}.$$
Thus,
$$\mathcal{D}_{u,v} = \bigcup_{0 \leq y \leq u} \mathcal{C}_{y,v}^+ = \bigcup_{x_0 \leq x \leq v} \mathcal{C}_{u,x}^-.$$

For $\beta = (\beta_1, \ldots, \beta_n) \in \mathbb{N}^n$, we set $\partial^\beta = \frac{\partial^{|\beta|}}{(\partial x)^{\beta_1}(\partial \theta^1)^{\beta_2} \ldots (\partial \theta^{n-1})^{\beta_n}}$ and recall $\partial_\mu = \frac{\partial}{\partial y^\mu}$ and (see (3.1)):

$$\begin{cases} \frac{\partial}{\partial x} = \frac{1}{2}\left(\frac{\partial}{\partial \tau} + \sum_{i=1}^{n} \frac{y^i}{\rho}\frac{\partial}{\partial y^i}\right) = \frac{1}{2}\left(\frac{\partial}{\partial \tau} + \frac{\partial}{\partial \rho}\right) \\ \frac{\partial}{\partial y} = \frac{1}{2}\left(\frac{\partial}{\partial \tau} - \sum_{i=1}^{n} \frac{y^i}{\rho}\frac{\partial}{\partial y^i}\right) = \frac{1}{2}\left(\frac{\partial}{\partial \tau} - \frac{\partial}{\partial \rho}\right) \end{cases} \quad (3.5)$$

Let $m \in \mathbb{N}$, $\alpha \in \mathbb{R}$ and $U$ a subset of $\mathbb{R}^{n+1}$. We will denote $H^m(U)$ the usual Sobolev space on $U$. Further for $U = \mathcal{C}^+$ or $\mathcal{C}_{u,v}^+$ we denote by

- $\mathscr{C}_0^\alpha(U)$ the set of continuous functions $f$ on $U$ for which the quantity $\|f\|_{\mathscr{C}_0^\alpha(U)} := \sup_{p \in U} |x|^{-\alpha}|f(p)|$ is finite,

- $\mathscr{C}_k^\alpha(U)$ the set of $k$−times continuously differentiable functions $f$ on $U$ such that the quantity $\|f\|_{\mathscr{C}_k^\alpha(U)} := \sum\limits_{0 \leq |\beta| \leq k} \||x|^{\beta_1} \partial^\beta f\|_{\mathscr{C}_0^\alpha(U)}$ is finite,

- $\mathscr{H}_k^\alpha(U)$ the space of those functions $f$ in $H^k_{loc}(U)$ for which the norm

$$\|f\|^2_{\mathscr{H}_k^\alpha(U)} := \sum_{0 \leq |\beta| \leq k} \int_U (|x|^{-\alpha+\beta_1} \partial^\beta f)^2 \frac{dx}{|x|} d\nu$$

is finite. Here $d\nu$ is a measure on $\mathscr{O}$ arising from a smooth Riemannian metric on $\mathbb{S}^{n-1}$.

## 3.3 Existence and uniqueness for a Goursat problem

Let $m \in \mathbb{N}$, $\alpha \in \mathbb{R}$. Let $\omega_0^-$, be a defined function on $\mathcal{C}^-$ and $\omega_0^+$, be defined on $\mathcal{C}^+$ such that

$$\omega_0^- \in H^{m+2}(\mathcal{C}^-); \quad \omega_0^+ \in \mathscr{H}_{m+1}^\alpha(\mathcal{C}^+) \quad \text{and} \quad \partial_x \omega_0^+ \in \mathscr{H}_{m+1}^\alpha(\mathcal{C}^+)^2, \tag{3.6}$$

and satisfying the compatibility condition

$$\omega_0^- = \omega_0^+ \quad \text{sur} \quad \mathcal{C}^+ \cap \mathcal{C}^- = \mathscr{O}. \tag{3.7}$$

The purpose of this section is to state and prove an existence and uniqueness theorem for the following characteristic Cauchy problem ($z = (y, x, \theta)$):

$$\begin{cases} \Box_{y,\eta}\omega = |x|^{-\frac{n+3}{2}} G\left(z, |x|^{\frac{n-1}{2}}\omega, |x|^{\frac{n-1}{2}}(\partial_y\omega, x\partial_x\omega, \partial_A\omega)\right) \text{ in } \mathcal{D} \\ \omega = \omega_0^+ \quad \text{on} \quad \mathcal{C}^+ \quad \text{and} \quad \omega = \omega_0^- \quad \text{on} \quad \mathcal{C}^- \end{cases}. \tag{3.8}$$

### 3.3.1 Hypothesis on the non linear term

In analogy with the procedure used in [6], we make the following assumption on the non linear source term $G$:

(H) We suppose that the function $G = G(z, p, q)$, is of $C^m$ class in all its variables and that the restriction $G(y, x)$ of $G$ on every slice $\{y = const\} \cap \{x = const\}$ has a *uniform zero of order* $r \geq 0$ at $p = q = 0$ in the sense that, for all $B > 0$ there exists a constant $\hat{C}(B)$ such that for $0 \leq j + \ell + i_1 + i_2 \leq \min(r, m)$ and $|(p, q)| \leq B$ one has:

$$\left\| \frac{|x|^{i_1} \partial^{j+\ell+i_1+i_2} G(y, x, \cdot, p, q)}{(\partial x)^{i_1}(\partial y)^{i_2} \partial p^j \partial q^\ell} \right\|_{C^{m-(j+\ell+i_1+i_2)}(\{y\}\times\{x\}\times\mathscr{O})} \leq \hat{C}(B) \|(p, q)\|^{r-j-\ell}. \tag{3.9}$$

**Remark 3.4** Note that if $G$ is of the form $G(z, p, q) = \sum\limits_{i+j=r}^{N} \Psi_{ij}(z) p^i q^j$ where the $\Psi_{ij}$'s are smooth functions then $G$ will satisfies hypothesis ($\mathcal{H}$).

**Remark 3.5** We point out for later use that this hypothesis implies that for all $\sigma \geq 0$, there exists a constant $C(\hat{C}, r, m, \sigma, B)$ such that for all $f \in H^m(\mathscr{O})$, with $\|f\|_{L^\infty(\mathscr{O})} \leq B$, we have

$$\|G(y, x, \cdot, |x|^\sigma f)\|_{H^m(\mathscr{O})} \leq C|x|^{r\sigma} \|f\|_{H^m(\mathscr{O})}. \tag{3.10}$$

---

[2] Actually, $\omega_0^+ \in \mathscr{H}_{m+1}^\alpha(\mathcal{C}^+)$ implies that $\partial_x \omega_0^+ \in \mathscr{H}_m^{\alpha-1/2}(\mathcal{C}^+)$, which will not be sufficient to obtain the control of some of our constants, we thus assume that $\partial_x \omega_0^+ \in \mathscr{H}_{m+1}^\alpha(\mathcal{C}^+)$.

### 3.3.2 First inequality

As a first step towards an existence theorem of the characteristic Cauchy problem (3.8), we prove now our first estimates. Let $\ell, \Lambda \in \mathbb{R}$, $\Lambda > 0$ and $\omega$ a sufficiently differentiable function defined on $\mathcal{Y}^+_{-\frac{1}{a},y}$. Set

$$L^\ell[\omega] = H(y,x)(\partial_x \omega + \partial_y \omega)\Box_{\eta,y}\omega \quad \text{with} \quad H(x,y) = (-x)^\ell e^{-\Lambda(y+x)}, \tag{3.11}$$

and

$$\nabla \omega = (\partial_x \omega, \partial_y \omega, \partial_\theta \omega), \; \nabla_x \omega = (\partial_x \omega, \partial_\theta \omega), \; \nabla_y \omega = (\partial_y \omega, \partial_\theta \omega)$$

where $\partial_\theta \omega = (\partial_{\theta^1}\omega, \ldots, \partial_{\theta^{n-1}}\omega)$. Assume that $c_0$ and $\bar{c}_0$ are two positive constants such that:

$$c_0 |X|^2 \leq X_x^2 + \sum_{A,B} \frac{h^{AB} X_A X_B}{2\rho^2} \leq \bar{c}_0 |X|^2. \tag{3.12}$$

We have the following:

**Proposition 3.6** *There exists a positive constant $c_1$ depending upon $h, c_0, \bar{c}_0$ and $n$ such that for all $\ell \geq 0$, $u \in [0, y_0]$, $v \in [x_0, 0[$, and for any function $\omega$ defined and at least of class $C^2$ on $\phi(\mathcal{Y}^+_{a,x}) = \mathcal{Y}^+_{-\frac{1}{a},y} \cap \mathcal{Y}^-_{0,y}$, we have:*

$$\int_{x_0}^v H(u,x) \|(\omega, \nabla_x \omega)(u,x)\|^2_{L^2(\mathcal{O})} dx + \int_0^u H(y,v) \|(\omega, \nabla_y \omega)(y,v)\|^2_{L^2(\mathcal{O})} dy \leq$$

$$\int_{x_0}^v H(0,x) \|(\omega, \nabla_x \omega)(0,x)\|^2_{L^2(\mathcal{O})} dx + \int_0^u H(y,x_0) \|(\omega, \nabla_y \omega)(y,x_0)\|^2_{L^2(\mathcal{O})} dy$$

$$+ (c_1(c_0, \bar{c}_0, n, h) - 2\Lambda) \int_0^u \int_{x_0}^v H(x,y) \|(\omega, \nabla \omega)(y,x)\|^2_{L^2(\mathcal{O})} dx dy$$

$$+ \frac{1}{c_0} \int_{\mathcal{D}_{u,v}} |L^\ell[\omega]| \, dy dx \, d\nu. \tag{3.13}$$

**Proof:** Recall that

$$\Box_{\eta,y}\omega = -4\partial_x \partial_y \omega + \frac{n-1}{\rho}(\partial_x - \partial_y)\omega + \frac{\Delta_{\mathbb{S}^{n-1}}\omega}{\rho^2}$$

$$= -4\partial_x \partial_y \omega + \frac{1}{\rho^2} h^{AB} \partial_A \partial_B \omega + \frac{n-1}{\rho}(\partial_x - \partial_y)\omega - \frac{1}{\rho^2}\Gamma^A \partial_A \omega,$$

where $\Gamma^A = h^{BC}\Gamma^A_{BC}$, the $\Gamma^A_{BC}$'s being the Christoffel symbols on the unit sphere $\mathbb{S}^{n-1}$ of $\mathbb{R}^n$. We point out the following trivial identities

$$-4H\partial_x \omega \partial_x \partial_y \omega = -2\partial_y \Big(H(\partial_x \omega)^2\Big) + (2\partial_y H)(\partial_x \omega)^2,$$

$$-4H\partial_y \omega \partial_x \partial_y \omega = -2\partial_x \Big(H(\partial_y \omega)^2\Big) + (2\partial_x H)(\partial_y \omega)^2,$$

$$\frac{1}{\rho^2} H h^{AB} \partial_x \omega \partial^2_{AB} \omega = \partial_A \left(\frac{1}{\rho^2} H h^{AB} \partial_x \omega \partial_B \omega\right) - \frac{1}{\rho^2} H \partial_A(h^{AB}) \partial_x \omega \partial_B \omega$$

$$- \frac{1}{\rho^2} H h^{AB} \partial_x \partial_A \omega \partial_B \omega.$$

The last term of the last identity can be written as

$$\frac{1}{\rho^2} H h^{AB} \partial_x \partial_A \omega \partial_B \omega = \frac{1}{2\rho^2} H h^{AB} (\partial_x \partial_A \omega \partial_B \omega + \partial_A \omega \partial_x \partial_B \omega)$$

$$= \frac{1}{2\rho^2} H h^{AB} \partial_x (\partial_A \omega \partial_B \omega)$$

$$= \partial_x \left(\frac{1}{2\rho^2} H h^{AB} \partial_A \omega \partial_B \omega\right) - \partial_x \left(\frac{1}{2\rho^2} H\right) h^{AB} \partial_A \omega \partial_B \omega.$$

It then follows that

$$\frac{1}{\rho^2} H h^{AB} \partial_x \omega \partial^2_{AB} \omega = \partial_A \left( \frac{1}{\rho^2} H h^{AB} \partial_x \omega \partial_B \omega \right) - \frac{1}{\rho^2} H \partial_A (h^{AB}) \partial_x \omega \partial_B \omega$$
$$- \partial_x \left( \frac{1}{2\rho^2} H h^{AB} \partial_A \omega \partial_B \omega \right) + \partial_x \left( \frac{1}{2\rho^2} H \right) h^{AB} \partial_A \omega \partial_B \omega \ .$$

Similar calculations lead to

$$\frac{1}{\rho^2} H h^{AB} \partial_y \omega \partial^2_{AB} \omega = \partial_A \left( \frac{1}{\rho^2} H h^{AB} \partial_y \omega \partial_B \omega \right) - \frac{1}{\rho^2} H \partial_A (h^{AB}) \partial_y \omega \partial_B \omega$$
$$- \partial_y \left( \frac{1}{2\rho^2} H h^{AB} \partial_A \omega \partial_B \omega \right) + \partial_y \left( \frac{1}{2\rho^2} H \right) h^{AB} \partial_A \omega \partial_B \omega \ .$$

We then obtain the following expression of $L^\ell[\omega]$ :

$$L^\ell[\omega] = -\partial_y \left( 2H(\partial_x \omega)^2 + \frac{1}{2\rho^2} H h^{AB} \partial_A \omega \partial_B \omega \right) - \partial_x \left( 2H(\partial_y \omega)^2 + \frac{1}{2\rho^2} H h^{AB} \partial_A \omega \partial_B \omega \right)$$
$$+ \partial_A \left( \frac{1}{\rho^2} H h^{AB} (\partial_B \omega)(\partial_x \omega + \partial_y \omega) \right) + \frac{n-1}{\rho} H \left( (\partial_x \omega)^2 - (\partial_y \omega)^2 \right)$$
$$+ 2 \partial_y H (\partial_x \omega)^2 + 2 \partial_x H (\partial_y \omega)^2 - \frac{1}{\rho^2} H \partial_A h^{AB} (\partial_B \omega)(\partial_x \omega + \partial_y \omega)$$
$$+ \frac{1}{2} h^{AB} \partial_A \omega \partial_B \omega (\partial_x + \partial_y)(H/\rho^2)) - \frac{1}{\rho^2} H \Gamma^A (\partial_x \omega + \partial_y \omega) \partial_A \omega \ .$$

Since $x = \tau + \rho$ and $y = \tau - \rho + \frac{1}{a}$ we have

$$\partial_x H = -\ell(-x)^{-1} H - \Lambda H, \ \partial_y H = -\Lambda H \text{ and } (\partial_x + \partial_y)(H/\rho^2) = \frac{1}{\rho^2} (-\ell(-x)^{-1} H - 2\Lambda H) \ . \quad (3.14)$$

Thus we have

$$L^\ell[\omega] = -\partial_y \left( 2H(\partial_x \omega)^2 + \frac{1}{2\rho^2} H h^{AB} \partial_A \omega \partial_B \omega \right) - \partial_x \left( 2H(\partial_y \omega)^2 + \frac{1}{2\rho^2} H h^{AB} \partial_A \omega \partial_B \omega \right)$$
$$+ \partial_A \left( \frac{1}{\rho^2} H h^{AB} (\partial_B \omega)(\partial_x \omega + \partial_y \omega) \right) + \frac{n-1}{\rho} H \left( (\partial_x \omega)^2 - (\partial_y \omega)^2 \right)$$
$$- \frac{1}{\rho^2} H \left( \partial_A h^{AB} + \Gamma^B \right) (\partial_B \omega)(\partial_x \omega + \partial_y \omega)$$
$$- 2\Lambda H \left( (\partial_x \omega)^2 + (\partial_y \omega)^2 + \frac{h^{AB}}{2\rho^2} \partial_A \omega \partial_B \omega \right)$$
$$- \ell(-x)^{-1} H \left( (\partial_y \omega)^2 + \frac{h^{AB}}{2\rho^2} \partial_A \omega \partial_B \omega \right) \ .$$

Integrating this identity on $\mathcal{D}_{u,v} = [0,u] \times [x_0, v] \times \mathcal{O}$ and using Stokes theorem, one is led to

$$\int_{\mathcal{D}_{u,v}} L^{\ell}[\omega] dy dx \, d\nu = -\int_{\partial \mathcal{D}_{u,v}} \left\{ 2(\partial_x \omega)^2 + \frac{h^{AB}}{2\rho^2} \partial_A \omega \partial_B \omega \right\} n_y H d\sigma$$
$$- \int_{\partial \mathcal{D}_{u,v}} \left\{ 2(\partial_y \omega)^2 + \frac{h^{AB}}{2\rho^2} \partial_A \omega \partial_B \omega \right\} n_x H d\sigma$$
$$+ \int_{\partial \mathcal{D}_{u,v}} \frac{h^{AB}}{\rho^2} (\partial_B \omega)(\partial_x \omega + \partial_y \omega) n_A H d\sigma$$
$$- \ell \int_{\mathcal{D}_{u,v}} (-x)^{-1} \left\{ (\partial_y \omega)^2 + \frac{h^{AB}}{2\rho^2} \partial_A \omega \partial_B \omega \right\} H dy dx \, d\nu$$
$$- 2\Lambda \int_{\mathcal{D}_{u,v}} \left\{ (\partial_x \omega)^2 + (\partial_y \omega)^2 + \frac{h^{AB}}{2\rho^2} \partial_A \omega \partial_B \omega \right\} H dy dx \, d\nu$$
$$+ (n-1) \int_{\mathcal{D}_{u,v}} \frac{1}{\rho} \left\{ (\partial_x \omega)^2 - (\partial_y \omega)^2 \right\} H dy dx \, d\nu$$
$$- \int_{\mathcal{D}_{u,v}} \frac{1}{\rho^2} \left( \partial_A h^{AB} + \Gamma^B \right) (\partial_B \omega)(\partial_x \omega + \partial_y \omega) H dy dx \, d\nu. \qquad (3.15)$$

In equation (3.15),

- $d\nu$ is the surface element on $\mathbb{S}^{n-1}$ defined by the induced metric on $\mathbb{S}^{n-1}$ by the Euclidean metric on $\mathbb{R}^n$,

- $n = n_y \partial_y + n_x \partial_x + \sum_{A=1}^{n} n_\theta \partial_\theta$ is the unit outward normal of the boundary $\partial \mathcal{D}_{u,v}$,

- and $d\sigma$ is the surface element on $\partial \mathcal{D}_{u,v}$ induced by the volume element $dy dx \, d\nu$.

The right-hand side of equation (3.15) is made of seven terms which will be labeled A, B, C, D, E, F and G where A is the terms of the first line, B those of the second line and so on.

**Remark 3.7** On the Riemannian manifold $\mathbb{R}^{n+1}$ endowed with the Euclidean metric, the family of vectors $\{\partial_\tau, \partial_\rho, \partial_\theta\}$ is an orthogonal frame and then we deduce that (note that $\partial \mathcal{D}_{u,v}$ is made of four pieces: $\partial \mathcal{D}_{u,v} = \mathcal{C}^+_{0,v} \cup \mathcal{C}^-_{u,x_0} \cup \mathcal{C}^+_{u,v} \cup \mathcal{C}^-_{u,v}$):

- on $\mathcal{C}^+_{0,v}$ the unit outward normal is $n = -\frac{1}{\sqrt{2}} \partial_y$, thus $n_y = -\frac{1}{\sqrt{2}}$, $n_x = 0$, $n_A = 0$, $A = 1, \ldots, n-1$;

- on $\mathcal{C}^+_{u,v}$ the outward unit normal is given by $n = \frac{1}{\sqrt{2}} \partial_y$, i.e. $n_y = \frac{1}{\sqrt{2}}$, $n_x = 0$, $n_A = 0$, $A = 1, \ldots, n-1$.

- on $\mathcal{C}^-_{u,x_0}$ the unit outward normal $n = -\frac{1}{\sqrt{2}} \partial_x$, i.e. $n_y = 0$, $n_x = -\frac{1}{\sqrt{2}}$, $n_A = 0$, $A = 1, \ldots, n-1$.

- on $\mathcal{C}^-_{u,v}$ we have $n = \frac{1}{\sqrt{2}} \partial_x$, thus $n_y = 0$, $n_x = \frac{1}{\sqrt{2}}$, $n_A = 0$, $A = 1, \ldots, n-1$.

In (3.15) we replace $\partial \mathcal{D}_{u,v}$ by $\mathcal{C}^+_{0,v} \cup \mathcal{C}^-_{u,x_0} \cup \mathcal{C}^+_{u,v} \cup \mathcal{C}^-_{u,v}$ and after using on each piece of $\partial \mathcal{D}_{u,v}$ the corresponding value of the outward unit normal, we find that:

$$A + B + C = \frac{1}{\sqrt{2}} \int_{\mathcal{C}^+_{0,v}} \left( 2(\partial_x \omega)^2 + \frac{h^{AB}}{2\rho^2} \partial_A \omega \partial_B \omega \right) H d\sigma - \frac{1}{\sqrt{2}} \int_{\mathcal{C}^+_{u,v}} \left( 2(\partial_x \omega)^2 + \frac{h^{AB}}{2\rho^2} \partial_A \omega \partial_B \omega \right) H d\sigma$$
$$+ \frac{1}{\sqrt{2}} \int_{\mathcal{C}^-_{u,x_0}} \left( 2(\partial_y \omega)^2 + \frac{h^{AB}}{2\rho^2} \partial_A \omega \partial_B \omega \right) H d\sigma - \frac{1}{\sqrt{2}} \int_{\mathcal{C}^-_{u,v}} \left( 2(\partial_y \omega)^2 + \frac{h^{AB}}{2\rho^2} \partial_A \omega \partial_B \omega \right) H d\sigma \, .$$

Identity (3.15) then takes the form:

$$\frac{1}{\sqrt{2}} \int_{\mathcal{C}_{u,v}^+} \left(2(\partial_x\omega)^2 + \frac{h^{AB}}{2\rho^2}\partial_A\omega\partial_B\omega\right)H d\sigma + \frac{1}{\sqrt{2}} \int_{\mathcal{C}_{u,v}^-} \left(2(\partial_y\omega)^2 + \frac{h^{AB}}{2\rho^2}\partial_A\omega\partial_B\omega\right)H d\sigma$$

$$= \frac{1}{\sqrt{2}} \int_{\mathcal{C}_{0,v}^+} \left(2(\partial_x\omega)^2 + \frac{h^{AB}}{2\rho^2}\partial_A\omega\partial_B\omega\right)H d\sigma + \frac{1}{\sqrt{2}} \int_{\mathcal{C}_{u,x_0}^-} \left(2(\partial_y\omega)^2 + \frac{h^{AB}}{2\rho^2}\partial_A\omega\partial_B\omega\right)H d\sigma$$

$$- \int_{\mathcal{D}_{u,v}} L^\ell[\omega] dy dx\, d\nu - \ell \int_{\mathcal{D}_{u,v}} (-x)^{-1}\left((\partial_y\omega)^2 + \frac{h^{AB}}{2\rho^2}\partial_A\omega\partial_B\omega\right)H dy dx\, d\nu$$

$$-2\Lambda \int_{\mathcal{D}_{u,v}} \left((\partial_x\omega)^2 + (\partial_y\omega)^2 + \frac{h^{AB}}{2\rho^2}\partial_A\omega\partial_B\omega\right)H dy dx\, d\nu + (n-1)\int_{\mathcal{D}_{u,v}} \frac{1}{\rho}\left((\partial_x\omega)^2 - (\partial_y\omega)^2\right)H dy dx\, d\nu$$

$$- \int_{\mathcal{D}_{u,v}} \frac{1}{\rho^2}\left(\partial_A h^{AB} + \Gamma^B\right)(\partial_B\omega)(\partial_x\omega + \partial_y\omega)H dy dx\, d\nu\,. \tag{3.16}$$

We then obtain the following estimate:

$$\int_{x_0}^v H(u,x)\|\nabla_x\omega(u,x)\|_{L^2(\mathscr{O})}^2 dx + \int_0^u H(y,v)\|\nabla_y\omega(y,v)\|_{L^2(\mathscr{O})}^2 dy$$

$$\leq \int_{x_0}^v H(0,x)\|\nabla_x\omega(0,x)\|_{L^2(\mathscr{O})}^2 dx + \int_0^u H(y,x_0)\|\nabla_y\omega(y,x_0)\|_{L^2(\mathscr{O})}^2 dy$$

$$+ \left(c(c_0, \bar{c}_0, n, \rho) - 2\Lambda\right)\int_0^u \int_{x_0}^v H(x,y)\|\nabla\omega(y,x)\|_{L^2(\mathscr{O})}^2 dx dy + \frac{1}{c_0}\int_{\mathcal{D}_{u,v}} \left|L^\ell[\omega]\right| dy dx\, d\nu\,. \tag{3.17}$$

On the other hand, we have:

$$\frac{1}{2}(\partial_x + \partial_y)(H\omega^2) = H\omega(\partial_x + \partial_y)\omega - \frac{\ell}{2}(-x)^{-1}H\omega^2 - \Lambda H\omega^2\,;$$

which implies that

$$(\partial_x + \partial_y)(H\omega^2) \leq 2H\omega(\partial_x + \partial_y)\omega - 2\Lambda H\omega^2$$
$$\leq \left((1-2\Lambda)\omega^2 + 2|\partial_x\omega|^2 + 2|\partial_y\omega|^2\right)H\,.$$

If we integrate once more on $\mathcal{D}_{u,v}$ then we obtain via Stokes formula the following inequality:

$$\int_{\partial\mathcal{D}_{u,v}} \omega^2 n_x H d\sigma + \int_{\partial\mathcal{D}_{u,v}} \omega^2 n_y H d\sigma \leq \int_{\mathcal{D}_{u,v}} \left((1-2\Lambda)\omega^2 + 2|\partial_x\omega|^2 + 2|\partial_y\omega|^2\right)H dy dx\, d\nu\,,$$

which is equivalent to

$$\int_{\mathcal{C}_{u,v}^+} \omega^2 H d\sigma + \int_{\mathcal{C}_{u,v}^-} \omega^2 H d\sigma \leq \int_{\mathcal{C}_{0,v}^+} \omega^2 H d\sigma + \int_{\mathcal{C}_{u,x_0}^-} \omega^2 H d\sigma$$
$$+ \sqrt{2}\int_{\mathcal{D}_{u,v}} \left((1-2\Lambda)\omega^2 + 2|\partial_x\omega|^2 + 2|\partial_y\omega|^2\right)H dy dx\, d\nu\,.$$

The estimate thus follows

$$\int_{x_0}^v H(u,x)\|\omega(u,x)\|_{L^2(\mathscr{O})}^2 dx + \int_0^u H(y,v)\|\omega(y,v)\|_{L^2(\mathscr{O})}^2 dy \leq$$

$$\int_{x_0}^v H(0,x)\|\omega(0,x)\|_{L^2(\mathscr{O})}^2 dx \int_0^{y_0} H(x_0,y)\|\omega(y,x_0)\|_{L^2(\mathscr{O})}^2 dy$$

$$+ \int_{\mathcal{D}_{u,v}} \left((1-2\Lambda)\omega^2 + 2|\partial_x\omega|^2 + 2|\partial_y\omega|^2\right)H dy dx\, d\nu\,. \tag{3.18}$$

Finally, adding side by side inequalities (3.17) and (3.18) leads to the stated inequality. $\square$

### 3.3.3 Iterative scheme

Our aim now, is to show that there exists a real number $u_* \in ]0, y_0]$ and a sequence of smooth functions $(\omega^k)_{k \in \mathbb{N}}$ which converges towards a solution $\omega$ of (3.8) on the set $\mathcal{D}_* := [0, u_*] \times [x_0, 0[ \times \mathcal{O}$. In order to use the $C^\infty$ results of [21], first, we need to approximate the data $\omega_0^+$ and $\omega_0^-$ with sequences of smooth functions $(\omega_0^{+,k})_{k \in \mathbb{N}}$ and $(\omega_0^{-,k})_{k \in \mathbb{N}}$ for which the compatibility condition

$$\omega_0^{+,k}(x_0, \theta) = \omega_0^{-,k}(0, \theta) \tag{3.19}$$

holds at every step of the iteration. These sequences are constructed as follows: denote by $(\bar{\omega}_0^{+,k})_{k \in \mathbb{N}}$ an arbitrary sequence of smooth functions which converges towards $\partial_x \omega_0^+$ in $\mathscr{H}_{m+1}^\alpha(\mathcal{C}^+)$ and by $(\omega_0^{-,k})_{k \in \mathbb{N}}$ an arbitrary sequence of smooth functions on $\mathcal{C}^-$ which converges to $\omega_0^-$ in the Sobolev space $H^{m+2}(\mathcal{C}^-)$. Then, for all $(x, \theta) \in [x_0, 0[ \times \mathcal{O}$ and $k \in \mathbb{N}$, set

$$\omega_0^{+,k}(x, \theta) = \omega_0^{-,k}(0, \theta) + \int_{x_0}^x \bar{\omega}_0^{+,k}(s, \theta) ds . \tag{3.20}$$

For later use we point out in the following Lemma some properties of the sequence $(\omega_0^{+,k})_{k \in \mathbb{N}}$.

**Lemma** 3.8 *Suppose that $-1 < \alpha \leq -1/2$. Then, the sequence $(\omega_0^{+,k})_{k \in \mathbb{N}}$ given by (3.20) satisfies the following:*

1. $\forall\ \theta \in \mathcal{O},\ \omega_0^{+,k}(x_0, \theta) = \omega_0^{-,k}(0, \theta)$ ;

2. $\omega_0^{+,k} \longrightarrow \omega_0^+$ in $\mathscr{H}_{m+1}^\alpha(\mathcal{C}^+)$ and $\partial_x \omega_0^{+,k} \longrightarrow \partial_x \omega_0^+$ in $\mathscr{H}_{m+1}^\alpha(\mathcal{C}^+)$ ;

3. $\displaystyle\sup_{k \in \mathbb{N},\ x \in [x_0, 0[} (-x)^{-\alpha} \|(\omega_0^{+,k}, \nabla_x \omega_0^{+,k})(x)\|_{H^{m-1}(\mathcal{O})} < \infty$ .

For the proof of this Lemma, the reader is refered to Section A page 20.

We denote by $\omega_0^k$ the continuous functions defined on $\mathcal{C}^+ \cup \mathcal{C}^-$ which coincide with $\omega_0^{+,k}$ on $\mathcal{C}^+$ and with $\omega_0^{-,k}$ on $\mathcal{C}^-$. The sequence $(\omega^k)_{k \in \mathbb{N}}$ is then constructed by induction:

- Set $\omega^0 = \omega_0$ where $\omega_0$ is a smooth function defined on $\mathcal{D}$ and which coincides with $\omega_0^0$ on $\mathcal{C}^+ \cup \mathcal{C}^-$.

- Then, let $\omega^{k+1}$ be defined by iteration as the solution of the linear characteristic Cauchy problem

$$\begin{cases} \Box_{y,\eta} \omega^{k+1} = |x|^{-\frac{n+3}{2}} G\left(z, |x|^{\frac{n-1}{2}}(\omega^k, \nabla \omega^k)\right) & \text{in } \mathcal{D} \\ \omega^{k+1} = \omega_0^{k+1} & \text{on } \mathcal{C}^+ \cup \mathcal{C}^- \end{cases} . \tag{3.21}$$

In order to prove existence of the sequence $(\omega^k)_{k \in \mathbb{N}}$, first we have to prove existence of the function $\omega_0$ used in the above iterative scheme. We define $\omega_0$ for any $(y, x, \theta) \in \mathcal{D}$ by setting

$$\omega_0(y, x, \theta) = \omega_0^{+,0}(x, \theta) + \omega_0^{-,0}(y, \theta) - \omega_0^{-,0}(0, \theta). \tag{3.22}$$

Next we have to justify existence of a smooth solution of (3.21). We quote Theorem 1 of [21]. Actually that reference gives a local solution on a neighborhood of the intersecting hypersurfaces, but in the case of the linear problem (3.21), we will obtain a global solution on $\mathcal{D}$.

### 3.3.4 Boundedness properties of $(\omega^k)_{k \in \mathbb{N}}$

Set

$$\begin{aligned} C_0 &= \sup_{k \in \mathbb{N},\ x \in [x_0, 0[} |x|^{-\alpha} \left\{ \|(\omega_0^{+,k}, \nabla_x \omega_0^{+,k})(x)\|_{W^{1,\infty}(\mathcal{O})} + \|\partial_y \omega^k(0, x)\|_{W^{1,\infty}(\mathcal{O})} \right\} \\ &+ \sup_{k \in \mathbb{N},\ y \in [0, y_0]} \left\{ (-x_0)^{-\alpha} \|\partial_y \omega_0^{-,k}(y)\|_{W^{1,\infty}(\mathcal{O})} \right\} \end{aligned}$$

We will show later that $C_0 < \infty$. We have the following

**Lemma 3.9** *Assume (3.6) and (3.7) with $-1 < \alpha \leq -1/2$ and $m > \frac{n+7}{2}$. If the source term $G$ satisfies hypothesis $(\mathcal{H})$ page 7 with a zero of order $r$ such that*

$$n \geq 1 + \frac{4}{r-1} - 2\alpha, \tag{3.23}$$

*then there exists a real number $u_* \in ]0, y_0]$ for which*

$$\sup_{k \in \mathbb{N},\, (y,x) \in [0, u_*] \times [x_0, 0[} |x|^{-\alpha} \|(\omega^k, \nabla \omega^k)(y, x)\|_{W^{1,\infty}(\mathscr{O})} < 2C_0. \tag{3.24}$$

The proof of this Lemma can be found in Section B page 22.

This lemma will be useful only if we prove that the constant $C_0$ is a finite quantity. We thus have to prove that the quantity $\sup_{k \in \mathbb{N},\, x \in \mathcal{C}^+} |x|^{-\alpha} \|\partial_y \omega^k(0, x)\|_{W^{1,\infty}(\mathscr{O})}$ is finite. This will be a consequence of the next Lemma. Set $\overline{H}(x) = e^{-\Lambda x} |x|^{-2\alpha}$ and

$$\hat{C}_0 := \sup_{k \in \mathbb{N},\, x \in [x_0, 0[} \overline{H}^{\frac{1}{2}}(x) \|(\omega_0^{+,k}, \nabla_x \omega_0^{+,k})(x)\|_{H^m(\mathscr{O})} < \infty \quad \text{see Lemma 3.8}$$

$$\widetilde{C}_0 := \sup_{k \in \mathbb{N}} \overline{H}^{\frac{1}{2}}(x_0) \|\partial_y \omega_0^{-,k}(0)\|_{H^{m-1}(\mathscr{O})} + \sup_{x \in [x_0, 0[} \overline{H}^{\frac{1}{2}}(x) \|\partial_y \omega^0(0, x)\|_{H^{m-1}(\mathscr{O})} < \infty.$$

Note that the constant $\check{C}(y_0, 0)$ is defined in Equation (B.25) page 30. Here and elsewhere we write $A \lesssim B$ if and only if there exists a constant $c > 0$ such that $A \leq cB$. We have the following:

**Lemma 3.10** *Under the hypotheses of Lemma 3.9, we have:*

$$\sup_{k \in \mathbb{N},\, x \in [x_0, 0[} \overline{H}^{\frac{1}{2}}(x) \|\partial_y \omega^k(0, x)\|_{H^{m-1}(\mathscr{O})} < 2(\hat{C}_0 + \widetilde{C}_0). \tag{3.25}$$

**Proof:** The proof will be carried out by induction on the integer $k$. By definition of the constants $\hat{C}_0$ and $\widetilde{C}_0$ the assumption is fulfilled when $k = 0$. Suppose that

$$\sup_{x \in [x_0, 0[} \overline{H}^{\frac{1}{2}}(x) \|\partial_y \omega^k)(0, x)\|_{H^{m-1}(\mathscr{O})} < 2(C_0 + \widetilde{C}_0).$$

We shall prove that this inequality remains true if we replace $k$ with $k+1$. If in Inequality B.14 page 28 we choose $\{y = 0\}$ we have (note that in (B.14) there is no $\epsilon$ in the second line but things can be arranged from (B.13) so as to get an $\epsilon$ there):

$$\begin{aligned}
H(0, x) \|\partial_y \omega^{k+1}(0, x)\|_{H^{m-1}(\mathscr{O})}^2 &\leq H(0, x_0) \|\partial_y \omega^{k+1}(0, x_0)\|_{H^{m-1}(\mathscr{O})}^2 \\
&\quad + \epsilon c_3(h, c_0, \bar{c}_0) \int_{x_0}^{x} H(0, s) \|\nabla_x \omega^{k+1}(0, s)\|_{H^m(\mathscr{O})}^2 ds \\
&\quad + \epsilon C(C_0) \int_{x_0}^{x} H(0, s) \|(\omega^k, \nabla \omega^k)(0, s)\|_{H^{m-1}(\mathscr{O})}^2 ds.
\end{aligned}$$

This can be rewritten as

$$\begin{aligned}
\overline{H}(x) \|\partial_y \omega^{k+1}(0, x)\|_{H^{m-1}(\mathscr{O})}^2 &\leq \overline{H}(x_0) \|\partial_y \omega_0^{-,k+1}(0)\|_{H^{m-1}(\mathscr{O})}^2 \\
&\quad + \epsilon c_3(h, c_0, \bar{c}_0) \int_{x_0}^{x} \overline{H}(s) \|\nabla_x \omega_0^{+,k+1}(s)\|_{H^m(\mathscr{O})}^2 ds \\
&\quad + \epsilon C(C_0) \int_{x_0}^{x} \overline{H}(s) \|(\omega_0^{+,k}, \nabla_x \omega_0^{+,k})(s)\|_{H^{m-1}(\mathscr{O})}^2 ds \\
&\quad + \epsilon C(C_0) \int_{x_0}^{x} \overline{H}(s) \|\partial_y \omega^k(0, s)\|_{H^{m-1}(\mathscr{O})}^2 ds,
\end{aligned}$$

which implies that:

$$\overline{H}(x)\|\partial_y\omega^{k+1}(0,x)\|^2_{H^{m-1}(\mathscr{O})} \leq \widetilde{C}_0^2 + \epsilon c_3(h,c_0,\bar{c}_0)|x_0|\hat{C}_0^2 + \epsilon C(C_0)|x_0|\hat{C}_0^2 + 4\epsilon C(C_0)|x_0|(\hat{C}_0+\widetilde{C}_0)^2$$
$$\leq 4(\hat{C}_0+\widetilde{C}_0)^2 \quad \text{since } \epsilon \text{ is sufficiently small}.$$

We then obtain

$$\sup_{x\in[x_0,0[} \overline{H}^{\frac{1}{2}}(x)\|\partial_y\omega^{k+1}(0,x)\|^2_{H^{m-1}(\mathscr{O})} \leq 2(\hat{C}_0+\widetilde{C}_0),$$

and the proof is complete. $\square$

**Lemma** 3.11 *Under the hypotheses of Lemma 3.9, there exists a constant $M_0 > 0$ such that:*

$$\sup_{k\in\mathbb{N},\ (y,x)\in[0,u_*]\times[x_0,0[} \|\omega^k(y,x)\|_{W^{1,\infty}(\mathscr{O})} < M_0. \tag{3.26}$$

**Proof:** let $\gamma \in \mathbb{N}^{n-1}$, such that $|\gamma| \in \{0,1\}$. By Lemma 3.9, for all $(y,x) \in [0,u_*] \times [x_0,0[$, we have:

$$|\partial_\theta^\gamma \omega^k(y,x)| \leq |\partial_\theta^\gamma \omega_0^{-,k}(y)| + \int_{x_0}^x |\partial_x\partial_\theta^\gamma \omega^k(y,s)|ds$$
$$\leq |\partial_\theta^\gamma \omega_0^{-,k}(y)| + \int_{x_0}^x |s|^\alpha |s|^{-\alpha}\|\partial_x\omega^k(s,x)\|_{C^1(\mathscr{O})}ds$$
$$\leq \underbrace{\sup_{|\gamma|\in\{0,1\},\ k\in\mathbb{N},\ y\in[0,u_*]} |\partial_\theta^\gamma \omega_0^{-,k}(y)| + 2C_0\int_{x_0}^0 |s|^\alpha ds}_{:=M_0}.$$

$\square$

We also have the following:

**Lemma** 3.12 *Under the hypotheses of the previous lemma, there exists two constants $M_1 > 0$ and $M_2 > 0$ such that:*

$$\sup_{k\in\mathbb{N},\ (y,x)\in[0,u_*]\times[x_0,0[} \|\omega^k(y,x)\|_{H^{m-1}(\mathscr{O})} < M_1 \tag{3.27}$$

*and*

$$\sup_{k\in\mathbb{N},\ (y,x)\in[0,u_*]\times[x_0,0[} |x|^{-\alpha}\|(\partial_x\omega^k,\partial_y\omega^k)(y,x)\|_{H^{m-1}(\mathscr{O})} < M_2. \tag{3.28}$$

The proof of this Lemma can be found in Section C page 34.

### 3.3.5 Convergence of the sequence $(\omega^k)_{k\in\mathbb{N}}$ and existence

Set $\delta\omega^k = \omega^{k+1} - \omega^k$ and $\delta\nabla\omega^k = \nabla\omega^{k+1} - \nabla\omega^k$. We have the following (recall that $\mathcal{D}_* = [0,u_*] \times [x_0,0[\times\mathscr{O})$:

**Lemma** 3.13 *Under the hypotheses of Lemma 3.11, even if it means to replace $(\omega^k)_{k\in\mathbb{N}}$ by one of its subsequences, there exists two real numbers $\sigma \in ]0,1[$ and $\varsigma > 0$ such that for all $k \geq 1$,*

$$\|(-x)^{-\alpha}(\delta\omega^k,\delta(\nabla\omega^k))\|_{L^2(\mathcal{D}_*)} \leq \frac{\varsigma}{2^k} + \sigma\|(-x)^{-\alpha}(\delta\omega^{k-1},\delta(\nabla\omega^{k-1}))\|_{L^2(\mathcal{D}_*)}. \tag{3.29}$$

The proof of this lemma is given in section D page 35.

Now we have all we need to show that the sequence $(\omega^k)_{k\in\mathbb{N}}$ converges towards a function $\omega$ of class $C^2$ on $\mathcal{D}_* = [0,u_*] \times [x_0,0[\times\mathscr{O}$ which is a solution of the characteristic initial value problem (3.8). We have the following consequence of the previous Lemmas.

**Corollary** 3.14 *There exists a continuous and bounded function $\omega$ on $\mathcal{D}_*$ such that $(\omega^k)_{k \in \mathbb{N}}$ converges to $\omega$ uniformly on any compact subset of $\mathcal{D}_*$.*

**Proof:** We point out the elementary fact: If $(U_n)_{n \in \mathbb{N}}$ is a sequence of positive real numbers satisfying $U_{n+1} \leq \alpha U_n + \frac{\beta}{2^n}$, then
$$U_n \leq \alpha^n U_0 + 2\beta \left( \frac{(1/2)^n - (\alpha)^n}{1 - 2\alpha} \right) .$$
Therefore, the series $\sum U_n$ will converge if $0 \leq \alpha < 1/2$. This remark and inequality (3.29) show that the series of functions $\sum e^{-\Lambda(y+x)/2} |x|^{-\alpha} (\delta \omega^k, \nabla \delta \omega^k)$ converges in the space $L^2(\mathcal{D}_*)$. Since the sequence of partial sums of this series write $S_k = e^{-\Lambda(y+x)/2}|x|^{-\alpha}\left((\omega^k, \nabla \omega^k) - (\omega^0, \nabla \omega^0)\right)$, the sequence $((\omega^k, \nabla \omega^k))_{k \in \mathbb{N}}$ converges to a function $(\omega_\varepsilon, \tilde{\omega}_\varepsilon)$ in the space $L^2(\mathcal{D}_{*,\varepsilon})$, with $\mathcal{D}_{*,\varepsilon} = [0, u_*] \times [x_0, -\varepsilon] \times \mathscr{O}$, for any $0 < \varepsilon < -x_0$. Note that the continuous embedding $L^2(\mathcal{D}_{*,\varepsilon}) \hookrightarrow \mathscr{D}'(\mathcal{D}_{*,\varepsilon})$ implies that $\tilde{\omega}_\varepsilon = \nabla \omega_\varepsilon$. We define $\omega$ by setting for any $(y, x, \theta) \in \mathcal{D}_*$, $\omega(y, x, \theta) = \omega_\varepsilon(y, x, \theta)$ if $(y, x, \theta) \in \mathcal{D}_{*,\varepsilon}$. First of all we need to prove that $\omega$ is a well defined function. Let $\varepsilon_1, \varepsilon_2 \in [0, -x_0[$ such that $\varepsilon_1 > \varepsilon_2$. Since $\mathcal{D}_{*,\varepsilon_1} \subset \mathcal{D}_{*,\varepsilon_2}$, $L^2(\mathcal{D}_{*,\varepsilon_2})$ embeds continuously in $L^2(\mathcal{D}_{*,\varepsilon_1})$ and then, the convergence of the sequence $(\omega^k, \nabla \omega^k)_{k \in \mathbb{N}}$ towards $(\omega_{\varepsilon_2}, \nabla \omega_{\varepsilon_2})$ in $L^2(\mathcal{D}_{*,\varepsilon_2})$, also holds in $L^2(\mathcal{D}_{*,\varepsilon_1})$. By uniqueness of limits of sequences in this space one is led to
$$(\omega_{\varepsilon_1}, \nabla \omega_{\varepsilon_1}) = (\omega_{\varepsilon_2}, \nabla \omega_{\varepsilon_2}) \text{ almost everywhere on } \mathcal{D}_{*,\varepsilon_1} . \tag{3.30}$$
Let $\varepsilon \in [0, -x_0[$. By Lemma 3.11, the sequence $(\omega^k)_{k \in \mathbb{N}}$ is uniformly bounded on $\mathcal{D}_*$ and therefore is uniformly bounded on $\mathcal{D}_{*,\varepsilon}$, by Lemma 3.9 page 13 there exists a contant $C = C(C_0, x_0, \varepsilon)$ such that
$$\sup_{k \in \mathbb{N}} \|\nabla \omega^k\|_{L^\infty(\mathcal{D}_{*,\varepsilon})} < C(C_0, x_0, \varepsilon) ,$$
thus the sequence $(\omega^k)_{k \in \mathbb{N}}$ is uniformly equicontinuous on $\mathcal{D}_{*,\varepsilon}$. Then, By Arzela-Ascoli theorem, there exists a subsequence $(\omega^{k_j})_{j \in \mathbb{N}}$ of the sequence $(\omega^k)_{k \in \mathbb{N}}$ which converges uniformly on the compact set $\mathcal{D}_{*,\varepsilon}$ to a continuous function $\omega'_\varepsilon$. The embedding $C^0(\mathcal{D}_{*,\varepsilon}) \hookrightarrow L^2(\mathcal{D}_{*,\varepsilon})$ proves that this convergence also holds in $L^2(\mathcal{D}_{*,\varepsilon})$ and by uniqueness of limits in $L^2(\mathcal{D}_{*,\varepsilon})$ we conclude that the equality in (3.30) holds everywhere and that:

- $\omega$ is a continuous function on $\mathcal{D}_*$ ,
- the sequence $(\omega^{k_j})_{j \in \mathbb{N}}$ uniformly converges to $\omega$ on any compact subset of $\mathcal{D}_* = [0, u_*] \times [x_0, 0[ \times \mathscr{O}$ .

It remains to prove that $\omega$ is bounded on $\mathcal{D}_*$. From the Sobolev embedding theorem (recall $m - 1 > \frac{n-1}{2} + 2$), by (3.27) we have:
$$\sup_{(y,x) \in [0, u_*] \times [x_0, 0[} \|\omega^k(y, x)\|_{C^2(\mathscr{O})} < M_1 \quad \forall k \in \mathbb{N} .$$
By taking the limits in this estimate we obtain that $\omega$ is a bounded function as well as its angular derivatives up to order two on $\mathcal{D}_*$. $\square$

**Lemma** 3.15 *$\forall s \in [0, m-2] \cap \mathbb{N}$, $(\omega^{k_j}(y, x))_{j \in \mathbb{N}}$ converges towards $\omega(y, x)$ in $H^s(\mathscr{O})$ uniformly in $(y, x)$ on $[0, u_*] \times [x_0, 0[$ and*
$$\omega \in \bigcap_{0 \leq s \leq m-2} C^0([0, u_*] \times [x_0, 0[; H^s(\mathscr{O})) .$$

**Proof:** By the previous corollary, for all $(y, x) \in [0, u_*] \times [x_0, 0[$, the sequence $(\omega^{k_j}(y, x))_{j \in \mathbb{N}}$ converges to $\omega(y, x)$ in $C^0(\mathscr{O})$ and since $C^0(\mathscr{O}) \hookrightarrow L^2(\mathscr{O})$, this convergence also holds in the space $L^2(\mathscr{O})$. On the other hand, by Lemma 3.12 the sequence $(\omega^{k_j}(y, x))_{j \in \mathbb{N}}$ is bounded in the Hilbert space $H^{m-1}(\mathscr{O})$, uniformly in $(y, x) \in [0, u_*] \times [x_0, 0[$. By weak compactness there exists a subsequence of $(\omega^{k_j}(y, x))_{j \in \mathbb{N}}$ denoted again by the same symbol which converges weakly to a function $\bar{\omega} \in H^{m-1}(\mathscr{O})$. This weak convergence also holds in $L^2(\mathscr{O})$, and by uniqueness of the weak limits, we obtain that $\omega(y, x) = \bar{\omega}(y, x) \in H^{m-1}(\mathscr{O})$. Now, we use the Interpolation Theorem (see for example [20], page 38) with $p = r = 2$, $s = m-1$, $u = \omega^{k_{j_1}} - \omega^{k_{j_2}}$, $j_1, j_2 \in \mathbb{N}$. We then obtain: $q = 2$ and that for all $i \in \{0, 1, \ldots, m-1\}$,

$$\begin{aligned}
\sum_{|\gamma|=i} \|\partial_\theta^\gamma (\omega^{k_{j_1}} - \omega^{k_{j_2}})\|_{L^2(\mathscr{O})}^2 &\leq c \left( \sum_{|\gamma|=m-1} \|\partial_\theta^\gamma (\omega^{k_{j_1}} - \omega^{k_{j_2}})\|_{L^2(\mathscr{O})} \right)^{\frac{2i}{m-1}} \|\omega^{k_{j_1}} - \omega^{k_{j_2}}\|_{L^2(\mathscr{O})}^{2(1-\frac{i}{m-1})} \\
&\leq C(M_1) \|\omega^{k_{j_1}} - \omega^{k_{j_2}}\|_{L^2(\mathscr{O})}^{2(1-\frac{i}{m-1})} .
\end{aligned}$$

This estimate implies that, if $s < m-1$, then the sequence $(\omega^{k_j}(y, x))_{j \in \mathbb{N}}$ is a Cauchy sequence in the Hilbert space $H^s(\mathscr{O})$ uniformly in $(y, x)$. $\square$

**Corollary** 3.16 *The following holds:*

- $\omega \in C^2([0, u_*] \times [x_0, 0[ \times \mathscr{O})$ ,

- $\omega$ *solves the characteristic initial value problem (3.8).*

The reader will find the proof of this corollary in Appendix E page 37.

### 3.3.6   Uniqueness and statement of the results

We are now going to show that the solution of (3.8) constructed in the previous section is the unique $C^2$ solution. Let $\omega_1, \omega_2$ be two functions of differentiability class $C^2$ on $\mathcal{D}_*$ both solutions of (3.8). Set $\delta\omega = \omega_2 - \omega_1$ and $\delta G(z) = G(z, |x|^{-\frac{n-1}{2}}(\omega_2, \nabla\omega_2)) - G(z, |x|^{-\frac{n-1}{2}}(\omega_1, \nabla\omega_1))$. It follows that $\delta\omega$ solves the characteristic initial value problem with vanishing [3] data

$$\begin{cases} \Box_{y,\eta} \delta w = (-x)^{-\frac{n+3}{2}} \delta G & \text{dans } \mathcal{D}_* \\ \delta\omega = 0 & \text{sur } \mathcal{C}^+ \cup \mathcal{C}^- \end{cases} . \tag{3.31}$$

We repeat the proof of Lemma 3.13 page 35 with instead $\delta\omega$ and obtain the following inequality which is the equivalent of (D.7) there:

$$\|H^{\frac{1}{2}}(y, x)(\delta\omega, \nabla\delta\omega)\|_{L^2(\mathcal{D}_*)}^2 \leq \sigma^2 \|H^{\frac{1}{2}}(y, x)(\delta\omega, \nabla\delta\omega)\|_{L^2(\mathcal{D}_*)}^2 .$$

This proves that $\omega_1 = \omega_2$ almost everywhere and since these functions are continuous functions they are equal everywhere. We have thus proved

**Theorem** 3.17 *Consider the characteristic initial value problem (3.8) on the subset $\mathcal{D} = [0, y_0] \times [x_0, 0[ \times \mathscr{O}$ of $\mathbb{R}_y^{n+1}$. Suppose that the initial data $\omega_0^+$ and $\omega_0^-$ satisfy (3.6) and (3.7) with $m \geq \frac{n+7}{2}$ and $-1 < \alpha \leq -1/2$. Moreover suppose that the nonlinear source term $G$ satisfies the nullity property $(\mathcal{H})$ page 7 with a uniform zero of order $r \geq 2$ being such that*

$$n \geq 1 + \frac{4}{r-1} - 2\alpha . \tag{3.32}$$

*Then there exists a positive real number $u_* \in ]0, y_0]$ and a unique function $\omega$ of differentiability class $C^2$ on $\mathcal{D}_* = [0, u_*] \times [x_0, 0[ \times \mathscr{O}$, solution of (3.8) with the following properties:*

- $\sup\limits_{z \in \mathcal{D}_*} |\omega(z)| < \infty$ ,

- $\sup\limits_{z \in \mathcal{D}_*} |x|^{-\alpha} |\nabla\omega(z)| < \infty$ ,

---

[3] We point out that the expression "vanishing data" here is not the same as in [10]. In that refenrence, this expression is used to say that the data and all their derivatives of all order in all directions vanish on the initial surface. In our case only $\delta\omega$ and it's interiour devivatives on $\mathcal{C}^+$ vanish on $\mathcal{C}^+$ and those of $\mathcal{C}^-$ vanish on $\mathcal{C}^-$; which is sufficient for our estimates.

- $\forall s \in [0, m-2] \cap \mathbb{N}, \ \omega \in C^0([0, u_*] \times [x_0, 0[; \ H^s(\mathscr{O}))$ .

A direct consequence of Theorem 3.17 is an existence and uniqueness result for the Cauchy problem (1.1) on the light cone. We want to solve this problem on a neighborhood of the entire cone. For this purpose, we need to make sure the data are such that, the problem at hand can be solved locally on a neighborhood $V_{0,x}$ of the tip of the initial cone and that the restriction of this local solution on any incoming cone intersecting this neighborhood is of $H^{m+2}$−regularity class (as in (3.6)). In order to obtain this local solution, we will use the result of [9] (see Théorème 2, page 47 of this reference). For any $\tau \in ]-\frac{1}{a}, 0]$ set (recall the vertex of the cone $\phi(\mathcal{C}_{a,x}^+)$ is $(-1/a, 0)$) :

$$Y^\tau = \{(y^\mu) \in \phi(\mathcal{Y}_{a,x}^+), \ -\frac{1}{a} \leq y^0 \leq \tau\} , \tag{3.33a}$$

$$C^\tau = \{(y^\mu) \in \phi(\mathcal{C}_{a,x}^+), \ -\frac{1}{a} \leq y^0 \leq \tau\} , \tag{3.33b}$$

and for any sufficiently small $\varepsilon > 0$ we set

$$\mathcal{C}^+(\varepsilon) = \phi\left(\mathcal{C}_{a,x}^+\right) \setminus C^{-\varepsilon} .$$

We have the following

**Theorem** 3.18 *Let $m \in \mathbb{N}$. Consider the characteristic initial value problem (1.1) on the light cone in the unbounded domain $\mathcal{Y}_{a,x}^+$ of $\mathbb{R}_x^{n+1}$. Assume that the source term $F$ is a smooth function of all its variables and that the initial data $\varphi$ are such that:*

- *there exists a real number $0 < \varepsilon_0 < \frac{1}{2a}$ such that $\hat{\varphi} = \left(\Omega^{-\frac{n-1}{2}} \varphi \circ \phi^{-1}\right)\Big|_{C^{-\varepsilon_0}}$ satisfies the last hypothesis $\mathcal{H}_{m+2}$ of [9, Théorème 2, p. 47],*

- *and $\forall \ \varepsilon \ \in ]0, \varepsilon_0]$,*

$$\Omega^{-\frac{n-1}{2}} \varphi \circ \phi^{-1}\Big|_{\mathcal{C}^+(\varepsilon)} \ ; \ \partial_x(\Omega^{-\frac{n-1}{2}} \varphi \circ \phi^{-1})\Big|_{\mathcal{C}^+(\varepsilon)} \in \mathscr{H}_{m+1}^\alpha(\mathcal{C}^+(\varepsilon)), \tag{3.34}$$

*where $m \geq \frac{n+7}{2}$ and $-1 < \alpha \leq -1/2$. Further assume that the function $\widetilde{F}$ of equation(3.3) has a uniform zero of order $\ell \geq 2$ satisfying*

$$n \geq 1 + \frac{4}{\ell - 1} - 2\alpha . \tag{3.35}$$

*Then, there exist three real numbers $a_0, C, R$ such that $a_0 > a, \ C, R > 0$ and a unique function $f$ of class $C^2$ on the future neighborhood $\mathscr{V}$ of $\mathcal{C}_{a,x}^+$ defined by $\mathscr{V} = \underset{a \leq b \leq a_0}{\cup} \mathcal{C}_{b,x}^+$ solution of (1.1) with the following decay property*

$$\forall (t, x^i) \in \mathscr{V}, \quad |f(t, x^i)| \leq C(t+r)^{-\frac{n-1}{2}}, \quad r = \sqrt{\sum_{i=1}^n (x^i)^2} > R$$

$$\forall (t, x^i) \in \mathscr{V}, \quad |\partial_t f(t, x^i)| \leq C(t+r)^{-\frac{n-1}{2}-\alpha}, \quad r > R$$

$$\forall (t, x^i) \in \mathscr{V}, \quad |\partial_r f(t, x^i)| \leq C(t+r)^{-\frac{n-1}{2}-\alpha}, \quad r > R .$$

**Proof:** Since $m + 2 > \frac{n}{2} + 1$, the last statement of hypotheses $\mathcal{H}_{m+2}$ of [9, Théorème 2, page 47] is fulfilled, assumes that $\hat{\varphi}$ can be decompose as $\hat{\varphi} = \bar{\varphi}\big|_{\phi(\mathcal{C}_{a,x}^+)} + \hat{\varphi}_1$ where $\bar{\varphi}$ is a polynomial function of degree $2(m+1)$ on $Y^\tau$, $\tau > 0$, and where $\varphi_1$ belongs to a weighted Sobolev space of differentiability class $2m + 3$ on $C^\tau$, the weight being choose so as to control the singularities at the tip of the cone. The results of this reference yield a neighborhood $V_{0,y}$ of the tip of the cone $\phi(\mathcal{C}_{a,x}^+)$ in $\phi(\mathcal{Y}_{a,x}^+)$ and a local solution $\hat{f}_0$ which restriction on any incoming cone intersecting $V_{0,y}$ is in the usual Sobolev space $H^{m+2}$. We then apply Theorem 3.17 to the Goursat problem (2.7)-(2.8) with $\omega_0^- = \hat{f}_0\big|_{\mathcal{C}^-}$ and $\omega_0^+ = \hat{\varphi}|_{\mathcal{C}^+}$ and obtain a bounded solution $\hat{f}$ of this problem on the future neighborhood

$$\mathcal{D}_* = [0, u_*] \times [x_0, 0[ \times \mathcal{O} = \bigcup_{0 \leq u \leq u_*} \mathcal{C}_{u,0}^+$$

of $\mathcal{C}^+ \cup \mathcal{C}^+$. Note that by uniqueness, $\hat{f}$ and $\hat{f}_0$ coincide on the intersection of $\mathcal{D}_*$ with $V_{0,y}$ the future neighborhood of the tip of the initial cone $\mathcal{C}_{x,a}^+$ on which we obtain from Dossa's results [9] the local solution $\hat{f}_0$. Therefore, there exists a constant $c > 0$ such that for all $(y^\mu) \in \mathcal{D}_*$, $|\hat{f}(y^\mu)| < c$. This estimate can be rewritten as $|\hat{f} \circ \phi(x^\mu)| < c$, for all $(x^\mu) \in \phi^{-1}(\mathcal{D}_*)$, i.e (see (2.3)); $|f(x^\mu)| < c |\Omega \circ \phi|^{-\frac{n-1}{2}}$. By the definition $\Omega$, (see (2.3)) we have

$$|\Omega| = |-\eta_{\alpha\beta} y^\alpha y^\beta| = |-(y^0)^2 + \rho^2| = \frac{1}{(t+r)(t-r)} \leq \frac{\tilde{c}}{t+r}.$$

This proves that for all $(t, x^i) \in \phi^{-1}(D_*)$, $|f(t, x^i)| \leq c(t+r)^{-\frac{n-1}{2}}$. Now according to some of our previous calculations, we have:

$$\frac{\partial}{\partial x^\mu} = -\Omega \frac{\partial}{\partial y^\mu} - 2 y_\mu y^\alpha \frac{\partial}{\partial y^\alpha}$$

$$= x(\frac{1}{a} - y)\frac{\partial}{\partial y^\mu} - 2 y_\mu \left( x \partial_x + (y - \frac{1}{a}) \partial_y \right).$$

This identity implies that (recall $t = x^0$ and $\tau = y^0 = -y_0$)

$$\frac{\partial}{\partial t} = x(\frac{1}{a} - y)\partial_\tau + 2\tau \left( x \partial_x + (y - \frac{1}{a}) \partial_y \right).$$

Using identities (3.1) and (3.5) leads to:

$$\frac{\partial}{\partial t} = x^2 \partial_x + \left( y - \frac{1}{a} \right)^2 \partial_y . \tag{3.37}$$

On the other hand, we have $\frac{\partial}{\partial r} = \frac{x^i}{r} \frac{\partial}{\partial x^i}$ and $\frac{x^i}{r} = -\frac{y^i}{\rho}$ thus,

$$\frac{\partial}{\partial r} = \Omega \frac{y^i}{\rho} \frac{\partial}{\partial y^i} + 2 y_i \frac{y^i}{\rho} \left( x \partial_x + (y - \frac{1}{a}) \partial_y \right)$$

$$= -x(\frac{1}{a} - y)\frac{\partial}{\partial \rho} + 2\rho \left( x \partial_x + (y - \frac{1}{a}) \partial_y \right).$$

Again from identities (3.1) and (3.5) we obtain

$$\frac{\partial}{\partial r} = x^2 \partial_x - \left( y - \frac{1}{a} \right)^2 \partial_y . \tag{3.38}$$

For all $(t, x^i)$ and $(\tau, y^i)$ such that $(\tau, y^i) = \phi(t, x^i)$, we have (recall $x \partial_x \Omega = \left( y - \frac{1}{a} \right) \partial_y \Omega = \Omega$)

$$\partial_t f(t, x^i) = x^2 \partial_x f \circ \phi^{-1}(\tau, y^i) + \left( y - \frac{1}{a} \right)^2 \partial_y f \circ \phi^{-1}(\tau, y^i)$$

$$= x^2 \partial_x \left( \Omega^{\frac{n-1}{2}} \hat{f}(\tau, y^i) \right) + \left( y - \frac{1}{a} \right)^2 \partial_y \left( \Omega^{\frac{n-1}{2}} \hat{f}(\tau, y^i) \right)$$

$$= \frac{n-1}{2} \left( x + y - \frac{1}{a} \right) \Omega^{\frac{n-1}{2}} \hat{f}(\tau, y^i) + x^2 \Omega^{\frac{n-1}{2}} \partial_x \hat{f}(\tau, y^i) + \left( y - \frac{1}{a} \right)^2 \Omega^{\frac{n-1}{2}} \partial_y \hat{f}(\tau, y^i) .$$

From Theorem 3.17 we know that for $r > R$,

$$|\hat{f}| \lesssim 1 \; , \;\; (-x)^{-\alpha}|\partial_x \hat{f}| \lesssim 1 \;\; \text{and} \;\; (-x)^{-\alpha}|\partial_y \hat{f}| \lesssim 1 \; .$$

Thus for all $(t, x^i)$ such that $r > R$, we have (recall $|\Omega| \lesssim \frac{1}{t+r}$)

$$\begin{aligned}|\partial_t f(t, x^i)| &\lesssim (t+r)^{-\frac{n-1}{2}} + (t+r)^{-\frac{n-1}{2}-2-\alpha} + (t+r)^{-\frac{n-1}{2}-\alpha}\\ &\lesssim (t+r)^{-\frac{n-1}{2}-\alpha}.\end{aligned}$$

The same holds for $|\partial_r f(t, x^i)|$. This proves that in general, the decay at infinity of the derivatives of the solution is not as fast as the decay of the solution itself and complete the proof. □

## 3.4 Application to wave maps

The aim of this section is to show that Theorem 3.18 applies to wave maps with source manifold the Minkowski space-time. Let $(\mathcal{N}, g)$ be a smooth Riemannian manifold with finite dimension $N$, we wish to find a map $f : (\mathbb{R}^{n+1}_x, \eta) \to (\mathcal{N}, g)$ solving the Cauchy problem for the wave map equation. As in [6], we will be interested in maps $f$ which have the property that $f$ approaches a constant map $f_0$ as $r$ tends to infinity along lightlike directions, $f_0(x^\mu) = p_0 \in \mathcal{N}$ for $x^\mu \in \mathbb{R}^{n+1}_x$. Introducing normal coordinate around $p_0$, we can write $f = f^a$, $a = 1, \ldots, N$, with the functions $f^a$ satisfying the following system of semi-linear partial differential equations

$$\Box_{\eta_x} f^a = F^a(f, \partial f) \; ; \tag{3.39}$$

with

$$F^a(f, \partial f) := -\eta^{\alpha\beta} \Gamma^a_{bc}(f) \frac{\partial f^b}{\partial x^\alpha} \frac{\partial f^c}{\partial x^\beta} \; ;$$

and where the $\Gamma^a_{bc}$'s are the Christoffel symbols of the metric $g$. Using as before the conformal transformation

$$\phi : \mathbb{R}^{n+1}_x \setminus \mathcal{C}_{0,x} \to \mathbb{R}^{n+1}_y \; \text{by} \; x^\alpha \mapsto y^\alpha := \frac{x^\alpha}{\eta_{\lambda\mu} x^\lambda x^\mu} \; , \; \alpha = 0, 1, \ldots, n \; .$$

and setting again $\Omega = -\eta_{\alpha\beta} y^\alpha y^\beta$; $\hat{f} = \Omega^{-\frac{n-1}{2}} f \circ \phi^{-1}$, (3.39) reads ( see (2.4), page 4):

$$\Box_{\eta_y} \hat{f}^a = \Omega^{-\frac{n+3}{2}} \widetilde{F}^a(\hat{f}, \partial_{y^\mu} \hat{f}), \tag{3.40}$$

with

$$\begin{aligned}\widetilde{F}^a(\hat{f}, \partial_{y^\mu} \hat{f}) &= -\Omega \Gamma^a_{bc}(\Omega^{\frac{n-1}{2}} \hat{f}) \Big\{ \Omega \eta^{\alpha\beta} (\Omega^{\frac{n-1}{2}} \partial_{y^\alpha} \hat{f}^b)(\Omega^{\frac{n-1}{2}} \partial_{y^\beta} \hat{f}^c) \\ &\quad - (1-n)^2 (\Omega^{\frac{n-1}{2}} \hat{f}^b)(\Omega^{\frac{n-1}{2}} \hat{f}^c) + 2(1-n)(\Omega^{\frac{n-1}{2}} \hat{f}^b) y^\mu (\Omega^{\frac{n-1}{2}} \partial_{y^\mu} \hat{f}^c) \Big\} \; .\end{aligned}$$

This expression shows that when transforming (3.40) with data on a null cone into a Goursat problem as in (2.7) we will instead have a pre-factor $\Omega^{-\frac{n+1}{2}}$. On the other hand, from the assumption on $f$, we know that $\widetilde{F}$ here has a *uniform zero* of order $r = 3$, thus in the case of wave maps, condition (3.23) reads:

$$n \geq 2 - 2\alpha \; .$$

We have proved the following:

**Theorem** 3.19 *Let $a > 0$, $n, m \in \mathbb{N}$, $n \geq 3$. Consider Equation (3.39) on the Minkowski space-time $\mathbb{R}^{n+1}_x$ with initial data given on the translated cone $\mathcal{C}^+_{a,x}$ and are such that:*

- *there exists a real number $0 < \varepsilon_0 < \frac{1}{2a}$ such that $\hat{\varphi} = \left( \Omega^{-\frac{n-1}{2}} f \circ \phi^{-1} \right)\Big|_{C^{\varepsilon_0}}$ satisfies the last hypothesis $\mathcal{H}_{m+2}$ of [9, Théorème 2, p. 47],*

- and $\forall \, \varepsilon \in \,]0, \varepsilon_0]$,

$$\Omega^{-\frac{n-1}{2}} f \circ \phi^{-1}|_{\mathcal{C}^+(\varepsilon)} \, ; \, \partial_x \left( \Omega^{-\frac{n-1}{2}} f \circ \phi^{-1}|_{\mathcal{C}^+(\varepsilon)} \right) \in \mathscr{H}^{-1/2}_{m+1}(\mathcal{C}^+(\varepsilon)) \, , \tag{3.41}$$

with $m \geq \frac{n+7}{2}$. Then, there exists three real numbers $a_0$, $C$, $R$ such that $a_0 > a$, $C, R > 0$ and a unique function $f$ of class $C^2$ on the future neighborhood $\mathscr{V}$ of $\mathcal{C}^+_{a,x}$ defined by $\mathscr{V} = \bigcup_{a \leq b \leq a_0} \mathcal{C}^+_{b,x}$ solution of (3.40) with the following decay property

$$\forall (t, x^i) \in \mathscr{V}, \quad |f(t, x^i)| \leq C(t+r)^{-\frac{n-1}{2}}, \quad r > R$$
$$\forall (t, x^i) \in \mathscr{V}, \quad |\partial_t f(t, x^i)| \leq C(t+r)^{-\frac{n-1}{2}-\alpha}, \quad r > R$$
$$\forall (t, x^i) \in \mathscr{V}, \quad |\partial_r f(t, x^i)| \leq C(t+r)^{-\frac{n-1}{2}-\alpha}, \quad r > R \, .$$

# A Proof of Proposition 3.8

The first statement is obvious. As far as the second statement is concerned, we write

$$\|\omega_0^{+,k} - \omega_0^+\|^2_{\mathscr{H}^\alpha_{m+1}(\mathcal{C}^+)} = \|\omega_0^{+,k} - \omega_0^+\|^2_{\mathscr{H}^\alpha_0(\mathcal{C}^+)} + \|x \partial_x (\omega_0^{+,k} - \omega_0^+)\|^2_{\mathscr{H}^\alpha_m(\mathcal{C}^+)} + \|\partial_A (\omega_0^{+,k} - \omega_0^+)\|^2_{\mathscr{H}^\alpha_m(\mathcal{C}^+)} \, .$$

We have

$$\|x \partial_x (\omega_0^{+,k} - \omega_0^+)\|^2_{\mathscr{H}^\alpha_m(\mathcal{C}^+)} = \|\bar{\omega}_0^{+,k} - \partial_x \omega_0^+\|^2_{\mathscr{H}^{\alpha-1/2}_m(\mathcal{C}^+)}$$
$$\leq c \|\bar{\omega}_0^{+,k} - \partial_x \omega_0^+\|^2_{\mathscr{H}^\alpha_m(\mathcal{C}^+)} \xrightarrow[k \to \infty]{} 0 \, . \tag{A.1}$$

On the other hand,

$$\omega_0^{+,k}(x, \theta) - \omega_0^+(x, \theta) = \omega_0^{-,k}(0, \theta) - \omega_0^+(x, \theta) + \int_{x_0}^x \bar{\omega}_0^{+,k}(s, \theta) ds$$
$$= \omega_0^{-,k}(0, \theta) - \omega_0^+(x_0, \theta) + \omega_0^+(x_0, \theta) - \omega_0^+(x, \theta) + \int_{x_0}^x \bar{\omega}_0^{+,k}(s, \theta) ds$$
$$= \omega_0^{-,k}(0, \theta) - \omega_0^-(0, \theta) + \int_{x_0}^x \left( \bar{\omega}_0^{+,k}(s, \theta) - \partial_x \omega_0^+(s, \theta) \right) ds \, ;$$

thus (recall $-1 < \alpha \leq -\frac{1}{2}$ implies $(-x)^{-2\alpha-1} \leq (-s)^{-2\alpha-1}$ for $s \leq x < 0$)

$$(-x)^{-2\alpha-1} |\omega_0^{+,k}(x, \theta) - \omega_0^+(x, \theta)|^2 \leq c(x_0) \left( |\omega_0^{-,k}(0, \theta) - \omega_0^-(0, \theta)|^2 + \int_{x_0}^x (-s)^{-2\alpha-1} |\bar{\omega}_0^{+,k} - \partial_x \omega_0^+|^2 (s, \theta) ds \right)$$
$$\leq c(x_0) \left( |\omega_0^{-,k}(0, \theta) - \omega_0^-(0, \theta)|^2 + \int_{x_0}^0 (-s)^{-2\alpha-1} |\bar{\omega}_0^{+,k} - \partial_x \omega_0^+|^2 (s, \theta) ds \right) \, .$$

Now integrating this inequality on $\mathcal{C}^+$ gives (the second inequality is obtained by trace theorem, see for example [14], Theorem 1 page 258):

$$\|\omega_0^{+,k} - \omega_0^+\|^2_{\mathscr{H}^\alpha_0(\mathcal{C}^+)} \leq c(x_0) \left( \|\omega_0^{-,k}(0) - \omega_0^-(0)\|^2_{L^2(\mathscr{O})} + \|\bar{\omega}_0^{+,k} - \partial_x \omega_0^+\|^2_{\mathscr{H}^\alpha_0(\mathcal{C}^+)} \right)$$
$$\leq c(x_0) \left( \|\omega_0^{-,k} - \omega_0^-\|^2_{H^1(\mathcal{C}^-)} + \|\bar{\omega}_0^{+,k} - \partial_x \omega_0^+\|^2_{\mathscr{H}^\alpha_0(\mathcal{C}^+)} \right) \xrightarrow[k \to \infty]{} 0 \, . \tag{A.2}$$

Now let $\beta \in \mathbb{N}^n$ such that $|\beta| \leq m$. Similarly to the previous calculations, we have

$$\partial^\beta \partial_A \left( \omega_0^{+,k}(x, \theta) - \omega_0^+(x, \theta) \right) = \partial^\beta \partial_A \left( \omega_0^{-,k}(0, \theta) - \omega_0^-(0, \theta) + \int_{x_0}^x \bar{\omega}_0^{+,k}(s, \theta) - \partial_x \omega_0^+(s, \theta) ds \right) \, .$$

If $\beta_1 = 0$ then,

$$\partial^\beta \partial_A \left(\omega_0^{+,k}(x,\theta) - \omega_0^+(x,\theta)\right) = \partial^\beta \partial_A \left(\omega_0^{-,k}(0,\theta) - \omega_0^-(0,\theta)\right) + \int_{x_0}^x \partial^\beta \partial_A \left(\bar{\omega}_0^{+,k}(s,\theta) - \partial_x \omega_0^+(s,\theta)\right) ds,$$

thus

$$(-x)^{-2\alpha-1} |\partial^\beta \partial_A \left(\omega_0^{+,k}(x,\theta) - \omega_0^+(x,\theta)\right)|^2$$
$$\leq c(x_0) \left(|\partial^\beta \partial_A \left(\omega_0^{-,k}(0,\theta) - \omega_0^-(0,\theta)\right)|^2 + \int_{x_0}^x (-s)^{-2\alpha-1}|\partial^\beta \partial_A(\bar{\omega}_0^{+,k} - \partial_x \omega_0^+)|^2(s,\theta) ds\right)$$
$$\leq c(x_0) \left(|\partial^\beta \partial_A \left(\omega_0^{-,k}(0,\theta) - \omega_0^-(0,\theta)\right)|^2 + \int_{x_0}^0 (-s)^{-2\alpha-1}|\partial^\beta \partial_A(\bar{\omega}_0^{+,k} - \partial_x \omega_0^+)|^2(s,\theta) ds\right).$$

Integrating on $\mathcal{C}^+$, we have:

$$\|\partial^\beta \partial_A(\omega_0^{+,k} - \omega_0^+)\|^2_{\mathscr{H}_0^\alpha(\mathcal{C}^+)} \leq c(x_0) \left(\|\omega_0^{-,k}(0) - \omega_0^-(0)\|^2_{H^{m+1}(\mathscr{O})} + \|\partial^\beta \partial_A(\bar{\omega}_0^{+,k} - \partial_x \omega_0^+)\|^2_{\mathscr{H}_0^\alpha(\mathcal{C}^+)}\right). \quad (A.3)$$

Suppose now $\beta_1 \geq 1$ and set $\tilde{\beta} = (\beta_1 - 1, \beta_2, \ldots, \beta_n)$. We have

$$\partial^\beta \partial_A \left(\omega_0^{+,k}(x,\theta) - \omega_0^+(x,\theta)\right) = \partial^{\tilde{\beta}} \partial_A \left(\bar{\omega}_0^{+,k}(s,\theta) - \partial_x \omega_0^+(s,\theta)\right),$$

thus,

$$\|(-x)^{\beta_1} \partial^\beta \partial_A \left(\omega_0^{+,k} - \omega_0^+\right)\|^2_{\mathscr{H}_0^\alpha} = \|(-x)^{\beta_1-1} \partial^{\tilde{\beta}} \partial_A \left(\bar{\omega}_0^{+,k} - \partial_x \omega_0^+\right)\|^2_{\mathscr{H}_0^{\alpha-1}}. \quad (A.4)$$

By (A.3) and (A.4), we have

$$\|\partial_A \left(\omega_0^{+,k}(x,\theta) - \omega_0^+\right)\|^2_{\mathscr{H}_m^\alpha} \leq c(x_0) \left(\|\omega_0^{-,k}(0) - \omega_0^-(0)\|^2_{H^{m+1}(\mathscr{O})} + \|\bar{\omega}_0^{+,k} - \partial_x \omega_0^+\|^2_{\mathscr{H}_{m+1}^\alpha}\right) \xrightarrow[k \to \infty]{} 0. \quad (A.5)$$

From (A.1), (A.2) and (A.5) it follows that $\omega_0^{+,k} \longrightarrow \omega_0^+$ in $\mathscr{H}_{m+1}^\alpha$. This proves that the sequence $(\omega_0^{+,k})_{k\in\mathbb{N}}$ is such that

- $\omega_0^{+,k} \longrightarrow \omega_0^+$ in $\mathscr{H}_{m+1}^\alpha(\mathcal{C}^+)$,
- $\partial_x \omega_0^{+,k} \longrightarrow \partial_x \omega_0^+$ in $\mathscr{H}_{m+1}^\alpha(\mathcal{C}^+)$.

Let now prove that the quantity $\sup_{k\in\mathbb{N},\ x\in[x_0,0[} (-x)^{-\alpha} \|\partial_x \omega_0^{+,k}(x)\|_{H^{m-1}(\mathscr{O})}$ is finite. We know that $[x_0,0[=\bigcup_{n\in\mathbb{N}^*} [\frac{x_0}{2^{n-1}}, \frac{x_0}{2^n}]$ and $s = \frac{2^n x}{x_0} \in [1,2]$ if and only if $x = \frac{sx_0}{2^n} \in [x_0, 0[$. For any function $f$ defined on $[x_0, 0[\times\mathscr{O}$, set $f_n(s,\theta) = f(x = \frac{sx_0}{2^n}, \theta)$. The $\mathscr{H}_m^\alpha(\mathcal{C}^+)$-norm of $f$ can be rewritten as (see Equation A.11 of [6]

$$\|f\|_{\mathscr{H}_m^\alpha(\mathcal{C}^+)} \approx (-x_0)^{-2\alpha} \sum_{n\geq 1} 2^{2n\alpha} \|f_n\|^2_{H^m([1,2]\times\mathscr{O})}. \quad (A.6)$$

Here we write $A \approx B$ if there exists constants $C_1, C_2 > 0$ such that $C_1 A \leq B \leq C_2 A$. We have the following

$$\sup_{x\in[x_0,0[} (-x)^{-2\alpha} \|f(x)\|^2_{H^{m-1}(\mathscr{O})} = \sup_{n\geq 1} \sup_{\frac{x_0}{2^{n-1}} \leq x \leq \frac{x_0}{2^n}} (-x)^{-2\alpha} \|f(x)\|^2_{H^{m-1}(\mathscr{O})}$$
$$= \sup_{n\geq 1} \sup_{s\in[1,2]} \left(\frac{sx_0}{2^n}\right)^{-2\alpha} \|f(\frac{sx_0}{2^n})\|^2_{H^{m-1}(\mathscr{O})}$$
$$= (-x_0)^{-2\alpha} \sup_{n\geq 1} \left\{2^{2n\alpha} \sup_{s\in[1,2]} (s)^{-2\alpha} \|f_n(s)\|^2_{H^{m-1}(\mathscr{O})}\right\}$$
$$\leq c(-x_0)^{-2\alpha} \sum_{n\geq 1} 2^{2n\alpha} \sup_{s\in[1,2]} \|f_n(s)\|^2_{H^{m-1}(\mathscr{O})}.$$

Now writing
$$\partial_\theta^\gamma f_n(s,\theta) = \partial_\theta^\gamma f_n(1,\theta) + \int_1^s \partial_s \partial_\theta^\gamma f_n(\tau,\theta)d\tau ,$$
implies that
$$|\partial_\theta^\gamma f_n(s,\theta)|^2 \leq |\partial_\theta^\gamma f_n(1,\theta)|^2 + c\int_1^s |\partial_s \partial_\theta^\gamma f_n(\tau,\theta)|^2 d\tau .$$

Integrating this estimate on $\mathcal{O}$ gives,
$$\begin{aligned}
\|f_n(s)\|^2_{H^{m-1}(\mathcal{O})} &\leq \|f_n(1)\|^2_{H^{m-1}(\mathcal{O})} + c\int_1^s \|\partial_s f_n(\tau)\|^2_{H^{m-1}(\mathcal{O})} d\tau \\
&\leq \|f_n(1)\|^2_{H^{m-1}(\mathcal{O})} + \|f_n\|^2_{H^m([1,2]\times\mathcal{O})} \\
&\leq c\|f_n\|^2_{H^m([1,2]\times\mathcal{O})} \quad \text{by trace theorem} .
\end{aligned}$$

Therefore,
$$\begin{aligned}
\sup_{x\in[x_0,0[} (-x)^{-2\alpha} \|f(x)\|^2_{H^{m-1}(\mathcal{O})} &\leq c(-x_0)^{-2\alpha} \sum_{n\geq 1} 2^{2n\alpha} \|f_n\|^2_{H^m([1,2]\times\mathcal{O})} \\
&\approx \|f_n\|_{\mathscr{H}^\alpha_m(\mathcal{C}^+)} \quad \text{(see (A.6))} .
\end{aligned}$$

Now choosing $f = (\omega_0^{+,k}, \nabla_x \omega_0^{+,k})$ in the previous estimate leads to
$$\sup_{x\in[x_0,0[} (-x)^{-2\alpha} \|(\omega_0^{+,k}, \nabla_x \omega_0^{+,k})(x)\|^2_{H^{m-1}(\mathcal{O})} \leq \|(\omega_0^{+,k}, \nabla_x \omega_0^{+,k})\|_{\mathscr{H}^\alpha_m(\mathcal{C}^+)} .$$

Since convergent sequences are bounded, we conclude that
$$\sup_{k\in\mathbb{N},\ x\in[x_0,0[} (-x)^{-2\alpha} \|(\omega_0^{+,k}, \nabla_x \omega_0^{+,k})(x)\|^2_{H^{m-1}(\mathcal{O})} < \infty .$$

This completes the proof of the Lemma 3.8. $\square$

# B  Proof of Lemma 3.9

The proof will be made by induction on the integer $k$. Let us show that the statement holds when $k=0$ i.e
$$\sup_{(y,x)\in[0,u_*]\times[x_0,0[} |x|^{-\alpha} \|(\omega^0, \nabla\omega^0)(y,x)\|_{W^{1,\infty}(\mathcal{O})} < 2C_0 .$$

From definition (3.22) of $\omega_0$, we have
$$\nabla\omega_0(y,x,\theta) = \left(\partial_y \omega_0^{-,0}(y,\theta), \partial_x \omega_0^{+,0}(x,\theta), \partial_\theta \omega_0^{+,0}(x,\theta) + \partial_\theta\left(\omega_0^{-,0}(y,\theta) - \omega^{-,0}(0,\theta)\right)\right)$$

thus,
$$\begin{aligned}
(-x)^{-\alpha} \|(\omega^0, \nabla\omega^0)(y,x)\|_{W^{1,\infty}(\mathcal{O})} &\leq (-x)^{-\alpha}\left(\|(\omega_0^{+,0}, \nabla_x\omega_0^{+,0})(x)\|_{W^{1,\infty}(\mathcal{O})} + \|\omega_0^{-,0}(y) - \omega_0^{-,0}(0)\|_{W^{1,\infty}(\mathcal{O})}\right) \\
&\quad +(-x)^{-\alpha}\left(\|\partial_y \omega_0^{-,0}(y)\|_{W^{1,\infty}(\mathcal{O})} + \|\partial_\theta(\omega_0^{-,0}(y) - \omega_0^{-,0}(0))\|_{W^{1,\infty}(\mathcal{O})}\right) \\
&\leq C_0 + (-x_0)^{-\alpha} \|\omega_0^{-,0}(y) - \omega_0^{-,0}(0)\|_{W^{1,\infty}(\mathcal{O})} \\
&\quad +(-x_0)^{-\alpha} \|\partial_\theta(\omega_0^{-,0}(y) - \omega_0^{-,0}(0))\|_{W^{1,\infty}(\mathcal{O})}
\end{aligned}$$

Now recall that $\omega_0^{-,0} \in C^\infty([0,y_0]\times\mathcal{O})$, thus
$$(-x_0)^{-\alpha}\|\omega_0^{-,0}(y) - \omega_0^{-,0}(0)\|_{W^{1,\infty}(\mathcal{O})} + (-x_0)^{-\alpha}\|\partial_\theta(\omega_0^{-,0}(y) - \omega_0^{-,0}(0))\|_{W^{1,\infty}(\mathcal{O})} \xrightarrow[y\to 0]{} 0 .$$

It then follows that there exists a real number $u_0 \in [0, y_0]$ such that, $\forall y \in [0, u_0]$,

$$(-x_0)^{-\alpha}\|(\omega_0^{-,0}(y) - \omega_0^{-,0}(0)\|_{W^{1,\infty}(\mathscr{O})} + (-x_0)^{-\alpha}\|\partial_\theta(\omega_0^{-,0}(y) - \omega_0^{-,0}(0))\|_{W^{1,\infty}(\mathscr{O})} < C_0 .$$

Therefore,

$$\sup_{(y,x)\in[0,u_0]\times[x_0,0[} |x|^{-\alpha}\|(\omega^0, \nabla\omega^0)(y,x)\|_{W^{1,\infty}(\mathscr{O})} < 2C_0 ,$$

and the property holds for $k = 0$. Note that the real $u_*$ will be determined later from the induction scheme and will be less than or equal to $u_0$. Let $j$ be an integer greater than or equal to 1, and suppose that for any integer $k \leq j$ the following holds:

$$\sup_{(y,x)\in[0,u_*]\times[x_0,0[} |x|^{-\alpha}\|(\omega^k, \nabla\omega^k)(y,x)\|_{W^{1,\infty}(\mathscr{O})} < 2C_0 . \tag{B.1}$$

We want to prove that (B.1) holds with $k = j+1$. Let $\gamma \in \mathbb{N}^{n-1}$ be such that $|\gamma| \leq m$. If in Proposition 3.6 page 8 we choose $\omega = \partial_\theta^\gamma \omega^{k+1}$ and $\ell = -2\alpha$ then we obtain the following inequality:

$$\int_{x_0}^v H(u,x)\|(\partial_\theta^\gamma \omega^{k+1}, \nabla_x \partial_\theta^\gamma \omega^{k+1})(u,x)\|_{L^2(\mathscr{O})}^2 dx + \int_0^u H(y,v)\|(\partial_\theta^\gamma \omega^{k+1}, \nabla_y \partial_\theta^\gamma \omega^{k+1})(y,v)\|_{L^2(\mathscr{O})}^2 dy \leq$$

$$\int_{x_0}^v H(0,x)\|(\partial_\theta^\gamma \omega^{k+1}, \nabla_x \partial_\theta^\gamma \omega^{k+1})(0,x)\|_{L^2(\mathscr{O})}^2 dx + \int_0^u H(y,x_0)\|(\partial_\theta^\gamma \omega^{k+1}, \nabla_y \partial_\theta^\gamma \omega^{k+1})(y,x_0)\|_{L^2(\mathscr{O})}^2 dy$$

$$+(c_1(c_0, \bar{c}_0, n, h) - 2\Lambda) \int_0^u \int_{x_0}^v H(x,y)\|(\partial_\theta^\gamma \omega^{k+1}, \nabla \partial_\theta^\gamma \omega^{k+1})(y,x)\|_{L^2(\mathscr{O})}^2 dxdy + \frac{1}{c_0}\int_{\mathcal{D}_{u,v}} |L^\ell[\partial_\theta^\gamma \omega^{k+1}]|\, dydx\, d\nu .$$

Summing up the above identities for all multi-indices $\gamma$ such that $|\gamma| \leq m$, one is led to:

$$\int_{x_0}^v H(u,x)\|(\omega^{k+1}, \nabla_x \omega^{k+1})(u,x)\|_{H^m(\mathscr{O})}^2 dx + \int_0^u H(y,v)\|(\omega^{k+1}, \nabla_y \omega^{k+1})(y,v)\|_{H^m(\mathscr{O})}^2 dy \leq$$

$$\int_{x_0}^v H(0,x)\|(\omega^{k+1}, \nabla_x \omega^{k+1})(0,x)\|_{H^m(\mathscr{O})}^2 dx + \int_0^u H(y,x_0)\|(\omega^{k+1}, \nabla_y \omega^{k+1})(y,x_0)\|_{H^m(\mathscr{O})}^2 dy$$

$$+(c(c_0, \bar{c}_0, n, \rho) - 2\Lambda) \int_0^u \int_{x_0}^v H(x,y)\|(\omega^{k+1}, \nabla \omega^{k+1})(y,x)\|_{H^m(\mathscr{O})}^2 dxdy$$

$$+\frac{1}{c_0}\sum_{|\gamma|\leq m}\int_{\mathcal{D}_{u,v}} |L^\ell[\partial_\theta^\gamma \omega^{k+1}]|\, dydx\, d\nu . \tag{B.2}$$

Let us control the terms with $L^\ell[\partial_\theta^\gamma \omega^{k+1}]$. In all the remaining of this section we will use the symbol $G^k(\ldots)$ to denote quantity $G\left(z, (-x)^{-\frac{n-1}{2}}(\omega^k, \nabla\omega^k)\right)$. We have:

$$\begin{aligned}
L^\ell[\partial_\theta^\gamma \omega^{k+1}] &= H(x,y)(\partial_x \partial_\theta^\gamma \omega^{k+1} + \partial_y \partial_\theta^\gamma \omega^{k+1})\Box_{\eta,y}\partial_\theta^\gamma \omega^{k+1} \\
&= H(x,y)(\partial_x \partial_\theta^\gamma \omega^{k+1} + \partial_y \partial_\theta^\gamma \omega^{k+1})\left(\partial_\theta^\gamma \Box_{\eta,y}\omega^{k+1} + [\Box_{\eta,y}, \partial_\theta^\gamma]\omega^{k+1}\right) \\
&= H(x,y)(\partial_x \partial_\theta^\gamma \omega^{k+1} + \partial_y \partial_\theta^\gamma \omega^{k+1})\left((-x)^{-\frac{n+3}{2}}\partial_\theta^\gamma G^k(\ldots) + [\Box_{\eta,y}, \partial_\theta^\gamma]\omega^{k+1}\right) \\
&=: A + B + C + D .
\end{aligned}$$

Here $C$ and $D$ are the terms with the commutators $[\Box_{\eta,y}, \partial_\theta^\gamma]$. We will use at many places the inequality $ab \leq a^2/(4\epsilon) + \epsilon b^2$. The term A is controlled as follows:

$$\begin{aligned}
A &= (-x)^{-\frac{n+3}{2}} H(x,y) \partial_\theta^\gamma \partial_x \omega^{k+1} \partial_\theta^\gamma G^k(\ldots) \\
&\leq c(\epsilon) H |\partial_\theta^\gamma \partial_x \omega^{k+1}|^2 + \epsilon H (-x)^{-(n+3)} |\partial_\theta^\gamma G^k(\ldots)|^2 ,
\end{aligned}$$

which implies

$$\sum_{|\gamma|\leq m}\int_{\mathcal{D}_{u,v}} A d\nu dx dy \leq c(\epsilon)\int_0^u\int_{x_0}^v H(x,y)\|\partial_x\omega^{k+1}(y,x)\|_{H^m(\mathscr{O})}^2 dxdy$$
$$+\epsilon\int_0^u\int_{x_0}^v (-x)^{-(n+3)} H(x,y)\|G^k(\ldots)\|_{H^m(\mathscr{O})}^2 dxdy\,.$$

Now, recall that from the induction hypothesis (B.1) we know that

$$\sup_{y,x}|x|^{-\alpha}\|(\omega^k,\nabla\omega^k)(y,x)\|_{L^\infty(\mathscr{O})} < 2C_0,$$

thus, one can use (3.10) to control the $H^m(\mathscr{O})$−norme of $G^k(\ldots)$:

$$\|G(\ldots)\|_{H^m(\mathscr{O})} = \|G\left(y,x,\theta,|x|^{\frac{n-1}{2}+\alpha}.|x|^{-\alpha}(w^k,\nabla_x w^k)\right)\|_{H^m(\mathscr{O})}$$
$$\leq C(C_0)|x|^{r\left(\frac{n-1}{2}+\alpha\right)}\||x|^{-\alpha}(\omega^k,\nabla\omega^k)(y,x)\|_{H^m(\mathscr{O})}\,.$$

Now, $-(n+3)+2r\left(\frac{n-1}{2}+\alpha\right)-2\alpha \geq 0$ if and only if $n \geq 1+\frac{4}{r-1}-2\alpha$. [4] We then obtain

$$\sum_{|\gamma|\leq m}\int_{\mathcal{D}_{u,v}} A d\nu dx dy \leq c(\epsilon)\int_0^u\int_{x_0}^v H(x,y)\|\partial_x\omega^{k+1}(y,x)\|_{H^m(\mathscr{O})}^2 dxdy$$
$$+\epsilon C(C_0)\int_0^u\int_{x_0}^v H(x,y)\|(\omega^k,\nabla\omega^k)(y,x)\|_{H^m(\mathscr{O})}^2 dxdy$$
$$\text{for}\quad n\geq 1+\frac{4}{r-1}-2\alpha\,. \qquad (B.3)$$

Similarly,

$$\sum_{|\gamma|\leq m}\int_{\mathcal{D}_{u,v}} B d\nu dx dy \leq c(\epsilon)\int_0^u\int_{x_0}^v H(x,y)\|\partial_y\omega^{k+1}(y,x)\|_{H^m(\mathscr{O})}^2 dxdy$$
$$+\epsilon C(C_0)\int_0^u\int_{x_0}^v H(x,y)\|(\omega^k,\nabla\omega^k)(y,x)\|_{H^m(\mathscr{O})}^2 dxdy$$
$$\text{for}\quad n\geq 1+\frac{4}{r-1}-2\alpha\,. \qquad (B.4)$$

As far as the terms $C$ and $D$ are concerned, we recall that the commutator reads:

$$[\Box_{\eta,y},\partial_\theta^\gamma]\omega^{k+1} = \frac{1}{\rho^2}h^{AB}\partial_A\partial_B\partial_\theta^\gamma\omega^{k+1} - \frac{1}{\rho^2}\Gamma^B\partial_B\partial_\theta^\gamma\omega^{k+1} - \frac{1}{\rho^2}\partial_\theta^\gamma\left(\frac{h^{AB}}{\rho^2}\partial_A\partial_B\omega^{k+1}\right) + \frac{1}{\rho^2}\partial_\theta^\gamma\left(\Gamma^B\partial_B\partial_\theta^\gamma\omega^{k+1}\right)$$
$$= -\sum_{\gamma_1\neq 0,\,\gamma_1+\gamma_2=\gamma} c(\gamma,\rho)\partial_\theta^{\gamma_1}h^{AB}\partial_\theta^{\gamma_2}\partial_A\partial_B\omega^{k+1} + \sum_{\gamma_1\neq 0,\,\gamma_1+\gamma_2=\gamma} c(\gamma,\rho)\partial_\theta^{\gamma_1}\Gamma^B\partial_\theta^{\gamma_2}\partial_B\omega^{k+1}\,;$$

whence, using inequality $ab \leq a^2+b^2$ one has:

$$\sum_{|\gamma|\leq m}\int_{\mathcal{D}_{u,v}} C d\nu dx dy \leq C(h,\rho)\int_0^u\int_{x_0}^v H(x,y)\|\nabla_x\omega^{k+1}(y,x)\|_{H^m(\mathscr{O})}^2 dxdy \qquad (B.5)$$

---

[4]The constraint $n\geq 1+\frac{4}{r-1}-2\alpha$ ensures that $-(n+3)+2r\left(\frac{n-1}{2}+\alpha\right)-2\alpha \geq 0$, and $(-x)^{-(n+3)+2r\left(\frac{n-1}{2}+\alpha\right)-2\alpha}$ is a bounded quantity in the range of coordinates we are concerned with.

and
$$\sum_{|\gamma|\leq m} \int_{\mathcal{D}_{u,v}} Dd\nu dxdy \leq C(h,\rho) \int_0^u \int_{x_0}^v H(x,y)\|\nabla_y \omega^{k+1}(y,x)\|_{H^m(\mathcal{O})}^2 dxdy \ . \tag{B.6}$$

Summing inequalities (B.3)-(B.6) gives:

$$\sum_{|\gamma|\leq m} \int_{\mathcal{D}_{u,v}} |L^\ell[\partial_\theta^\gamma \omega^{k+1}]| \, dydx\, d\nu \leq c(\epsilon) \int_0^u \int_{x_0}^v H(x,y)\|\nabla \omega^{k+1}(y,x)\|_{H^m(\mathcal{O})}^2 dxdy$$
$$+\epsilon C(C_0) \int_0^u \int_{x_0}^v H(x,y)\|(\omega^k, \nabla\omega^k)(y,x)\|_{H^m(\mathcal{O})}^2 dxdy$$
$$\text{for} \quad n \geq 1 + \frac{4}{r-1} - 2\alpha \ .$$

We can then rewrite (B.2) as:

$$\int_{x_0}^v H(u,x)\|(\omega^{k+1}, \nabla_x \omega^{k+1})(u,x)\|_{H^m(\mathcal{O})}^2 dx + \int_0^u H(y,v)\|(\omega^{k+1}, \nabla_y \omega^{k+1})(y,v)\|_{H^m(\mathcal{O})}^2 dy \leq$$
$$\int_{x_0}^v H(0,x)\|(\omega^{k+1}, \nabla_x \omega^{k+1})(0,x)\|_{H^m(\mathcal{O})}^2 dx + \int_0^u H(y,x_0)\|(\omega^{k+1}, \nabla_y \omega^{k+1})(y,x_0)\|_{H^m(\mathcal{O})}^2 dy$$
$$+(c(c_0,\bar{c}_0,n,\rho) + c(c_0,\epsilon) - 2\Lambda) \int_0^u \int_{x_0}^v H(x,y)\|(\omega^{k+1}, \nabla\omega^{k+1})(y,x)\|_{H^m(\mathcal{O})}^2 dxdy$$
$$+\epsilon C(C_0) \int_0^u \int_{x_0}^v H(x,y)\|(\omega^k, \nabla\omega^k)(y,x)\|_{H^m(\mathcal{O})}^2 dxdy$$
$$\text{for} \quad n \geq 1 + \frac{4}{r-1} - 2\alpha \ .$$

All the derivatives appearing in the first term of the second line of the above equation are interior derivatives to the hypersurface $\{y=0\}$ and those of the second term are interior to $\{x=x_0\}$, therefore we can rewrite this last estimate using the initial data of the Cauchy problem (3.21):

$$\int_{x_0}^v H(u,x)\|(\omega^{k+1}, \nabla_x \omega^{k+1})(u,x)\|_{H^m(\mathcal{O})}^2 dx + \int_0^u H(y,v)\|(\omega^{k+1}, \nabla_y \omega^{k+1})(y,v)\|_{H^m(\mathcal{O})}^2 dy$$
$$\leq \int_{x_0}^v H(0,x)\|(\omega_0^{+,k+1}, \nabla_x \omega_0^{+,k+1})(x)\|_{H^m(\mathcal{O})}^2 dx + \int_0^u H(y,x_0)\|(\omega_0^{-,k+1}, \nabla_y \omega_0^{-,k+1})(y)\|_{H^m(\mathcal{O})}^2 dy$$
$$+(c(c_0,\bar{c}_0,n,\rho) + c(c_0,\epsilon) - 2\Lambda) \int_0^u \int_{x_0}^v H(x,y)\|(\omega^{k+1}, \nabla\omega^{k+1})(y,x)\|_{H^m(\mathcal{O})}^2 dxdy$$
$$+\epsilon C(C_0) \int_0^u \int_{x_0}^v H(x,y)\|(\omega^k, \nabla\omega^k)(y,x)\|_{H^m(\mathcal{O})}^2 dxdy$$
$$\text{for} \quad n \geq 1 + \frac{4}{r-1} - 2\alpha \ .$$

We choose now in the above estimate the parameter $\Lambda = \Lambda_0$ large enough so that $c(c_0,\bar{c}_0,n,\rho) + c(c_0,\epsilon) - 2\Lambda_0 < 0$ and we have prove the following

**Lemma** B.1 *Suppose $n \geq 1 + \frac{4}{r-1} - 2\alpha$. For all $\epsilon \in ]0,1]$ there exists a constant $\Lambda_0 = \Lambda_0(c, c_0, \bar{c}_0, \epsilon, h) > 0$* such that $\forall (u,v) \in [0, y_0] \times [x_0, 0[$ and $k \leq j$, we have:

$$\int_{x_0}^{v} H(u,x)\|(\omega^{k+1}, \nabla_x \omega^{k+1})(u,x)\|_{H^m(\mathscr{O})}^2 dx + \int_{0}^{u} H(y,v)\|(\omega^{k+1}, \nabla_y \omega^{k+1})(y,v)\|_{H^m(\mathscr{O})}^2 dy \leq$$

$$\int_{x_0}^{v} H(0,x)\|(\omega_0^{+,k+1}, \nabla_x \omega_0^{+,k+1})(x)\|_{H^m(\mathscr{O})}^2 dx + \int_{0}^{u} H(y,x_0)\|(\omega_0^{-,k+1}, \nabla_y \omega_0^{-,k+1})(y)\|_{H^m(\mathscr{O})}^2 dy$$

$$+\epsilon C(C_0) \int_{0}^{u} \int_{x_0}^{v} H(x,y)\|(\omega^k, \nabla \omega^k)(y,x)\|_{H^m(\mathscr{O})}^2 dx dy \ . \tag{B.7}$$

One would like to get rid of the dependence of $k$ in the right-hand side of the above estimate. We proceed as follows. For any $(u,v) \in [0, y_0] \times [x_0, 0]$, set

$$\begin{aligned}
\hat{C}(u,v) &= |x_0|^{-2\alpha} e^{-\Lambda_0 x_0} \sup_{k \in \mathbb{N}} \int_{0}^{u} e^{-\Lambda_0 y} \|(\omega_0^{-,k}, \nabla_y \omega_0^{-,k})(y)\|_{H^m(\mathscr{O})}^2 dy \\
&\quad + \sup_{k \in \mathbb{N}} \int_{x_0}^{v} |x|^{-2\alpha} e^{-\Lambda_0 x} \|(\omega_0^{+,k}, \nabla_x \omega_0^{+,k})(x)\|_{H^m(\mathscr{O})}^2 dx \\
&\quad + \frac{1}{2(y_0 + |x_0|)} \int_{0}^{u} \int_{x_0}^{v} H(y,x) \|(\omega^0, \nabla \omega^0)\|_{H^m(\mathscr{O})}^2 \ . 
\end{aligned} \tag{B.8}$$

We notice that $\forall (u,v) \in [0, y_0] \times [x_0, 0]$, the quantity $\hat{C}(u,v)$ is finite. Indeed we have

$$\begin{aligned}
\hat{C}(u,v) &\leq c(x_0, y_0, \Lambda_0) \left( \sup_{k \in \mathbb{N}} \|\omega_0^{-,k}\|_{H^{m+1}(\mathcal{C}^-)}^2 + \sup_{k \in \mathbb{N}} \|(\omega_0^{+,k}, \nabla_x \omega_0^{+,k})\|_{\mathscr{H}_m^\alpha(\mathcal{C}^+)}^2 \right) \\
&\quad + c(x_0, y_0, \Lambda_0) \int_{0}^{u} \int_{x_0}^{v} |x|^{-2\alpha} \|(\omega^0, \nabla \omega^0)\|_{H^m(\mathscr{O})}^2 \ .
\end{aligned}$$

The two terms in the first line of this estimate are bounded because convergent sequences (see Lemma 3.8) are bounded. From (3.22) we have

$$\int_{0}^{u} \int_{x_0}^{v} |x|^{-2\alpha} \|(\omega^0, \nabla \omega^0)\|_{H^m(\mathscr{O})}^2 \leq C \Big( u \|(\omega_0^{0,+}, \nabla_x \omega_0^{0,+})\|_{\mathscr{H}_m^\alpha(\mathcal{C}^+)}^2 + (v - x_0) \|(\omega_0^{0,-}, \nabla_y \omega_0^{0,-})\|_{H^m(\mathcal{C}^-)}^2$$

$$+ u(v - x_0) \|\omega_0^{-,0}(0)\|_{H^{m+1}(\mathscr{O})}^2 \Big) < \infty \ .$$

This proves that (B.8) defines a finite quantity. Now by the definition of this constant, (B.7) implies:

$$\int_{x_0}^{v} H(u,x) \|(\omega^{k+1}, \nabla_x \omega^{k+1})(u,x)\|_{H^m(\mathscr{O})}^2 dx \leq \hat{C}(u,v)$$

$$+ \epsilon C(C_0, c_0) \int_{0}^{u} \int_{x_0}^{v} H(x,y) \|(\omega^k, \nabla \omega^k)(y,x)\|_{H^m(\mathscr{O})}^2 dx dy \ . \tag{B.9}$$

Suppose that for all $u, v \in [0, y_0] \times [x_0, 0]$

$$\int_{0}^{u} \int_{x_0}^{v} H(x,y) \|(\omega^k, \nabla \omega^k)(y,x)\|_{H^m(\mathscr{O})}^2 dx dy \leq 2\hat{C}(y_0, 0)(y_0 + |x_0|) \ . \tag{B.10}$$

After integration with respect to $y$ on $[0, \tilde{u}]$ for any $\tilde{u} \in [0, y_0]$, inequality (B.9) gives

$$\int_{0}^{\tilde{u}} \int_{x_0}^{v} H(u,x) \|(\omega^{k+1}, \nabla_x \omega^{k+1})(u,x)\|_{H^m(\mathscr{O})}^2 dx du \leq \hat{C}(y_0, 0) y_0 \Big( 1 + 2\epsilon C(C_0, c_0)(y_0 + |x_0|) \Big)$$

$$\leq 2\hat{C}(y_0, 0) y_0 \ ,$$

if $\epsilon$ small enough. Using now this inequality in (B.7) leads to:

$$\int_0^u H(y,v)\|(\omega^{k+1}, \nabla_y \omega^{k+1})(y,v)\|^2_{H^m(\mathcal{O})} dy \leq \hat{C}(y_0, 0) + 2\epsilon C(C_0, c_0)\hat{C}(y_0, 0)(y_0 + |x_0|) \ .$$

Again by integration, for any $\tilde{v} \in [x_0, 0]$ we have

$$\int_0^u \int_{x_0}^{\tilde{v}} H(y,v)\|(\omega^{k+1}, \nabla_y \omega^{k+1})(y,v)\|^2_{H^m(\mathcal{O})} dy dv \leq \hat{C}(y_0, 0)\Big(1 + 2\epsilon C(C_0, c_0)(y_0 + |x_0|)\Big)|x_0|$$

$$\leq 2\hat{C}(y_0, 0)|x_0| \quad \text{if } \epsilon \text{ is small enough} \ .$$

Therefore, assuming that (B.10) is true, we have proved that

$$\int_0^u \int_{x_0}^v H(x,y)\|(\omega^{k+1}, \nabla \omega^{k+1})(y,x)\|^2_{H^m(\mathcal{O})} dx dy \leq 2\hat{C}(y_0, 0)(y_0 + |x_0|) \ .$$

Considering the definition of the constant $\hat{C}(u,v)$ (see (B.8)), we then obtain that (B.10) holds for $k=0$, and one can conclude that for any $k \in \mathbb{N}$ inequality (B.10) is satisfied. We have then proved

**Lemma** B.2 *Suppose that the constant $\hat{C}$ is defined by (B.8). One can choose $\epsilon = \epsilon_0(c_0, \bar{c}_0, x_0, y_0, C_0, h)$ such that*

$$\sup_{k \in \mathbb{N}, \ (u,v) \in [0,y_0] \times [x_0, 0[} \int_0^u \int_{x_0}^v H(x,y)\|(\omega^k, \nabla \omega^k)(y,x)\|^2_{H^m(\mathcal{O})} dx dy \leq 2\hat{C}(y_0, 0)(y_0 + |x_0|) \ , \tag{B.11}$$

*and for any $\Lambda \geq \Lambda_0$,*

$$\int_{x_0}^v H(u,x)\|(\omega^{k+1}, \nabla_x \omega^{k+1})(u,x)\|^2_{H^m(\mathcal{O})} dx + \int_0^u H(y,v)\|(\omega^{k+1}, \nabla_y \omega^{k+1})(y,v)\|^2_{H^m(\mathcal{O})} dy \leq 2\hat{C}(y_0, 0) \ . \tag{B.12}$$

**Remark** B.3 Note that as we assume that the induction hypothesis holds for any $k \leq j$, inequality (B.12) hold for any $k \leq j$.

Let $\gamma \in \mathbb{N}^{n-1}$, such that $|\gamma| \leq m-1$. To proceed further we apply $\partial_\theta^\gamma$ to the differential equation satisfies by $\omega^{k+1}$ and then multiply the differentiated equation by $H\partial_\theta^\gamma \partial_y \omega^{k+1}$. As in the proof of the Proposition 3.6, we obtain:

$$\begin{aligned}
\partial_x \left(H(\partial_y \partial_\theta^\gamma \omega^{k+1})^2\right) &= (\partial_x H - H\frac{n-1}{2\rho})(\partial_y \partial_\theta^\gamma \omega^{k+1})^2 + H\frac{n-1}{2\rho} \partial_y \partial_\theta^\gamma \omega^{k+1} \partial_x \partial_\theta^\gamma \omega^{k+1} \\
&\quad + \partial_y \partial_\theta^\gamma \omega^{k+1} \sum_{\gamma_1 + \gamma_2 = \gamma} H \frac{\partial_\theta^{\gamma_1} h^{AB}}{2\rho^2} \partial^{\gamma_2} \partial_A \partial_B \omega^{k+1} - \partial_y \partial_\theta^\gamma \omega^{k+1} \sum_{\gamma_1 + \gamma_2 = \gamma} H \frac{\partial_\theta^{\gamma_1} \Gamma^B}{2\rho^2} \partial_\theta^{\gamma_2} \partial_B \omega^{k+1} \\
&\quad - \frac{1}{2}|x|^{-\frac{n+3}{2}} H \partial_y \partial_\theta^\gamma \omega^{k+1} \partial_\theta^\gamma G^k(\dots) \ .
\end{aligned} \tag{B.13}$$

We integrate this identity on the set $\{y\} \times [x_0, x] \times \mathcal{O}$. From Stokes' theorem we have for $n \geq 1 + \frac{4}{r-1} - 2\alpha$:

$$\begin{aligned}
H(y,x)\|\partial_y \omega^{k+1}(y,x)\|^2_{H^{m-1}(\mathcal{O})} &\leq H(y,x_0)\|\partial_y \omega^{k+1}(y,x_0)\|^2_{H^{m-1}(\mathcal{O})} \\
&\quad + (c_2(h, c_0, \bar{c}_0) + c(\epsilon) - \Lambda) \int_{x_0}^x H(y,s)\|\partial_y \omega^{k+1}(y,s)\|^2_{H^{m-1}(\mathcal{O})} ds \\
&\quad + c_3(h, c_0, \bar{c}_0) \int_{x_0}^x H(y,x)\|\nabla_x \omega^{k+1}(y,s)\|^2_{H^m(\mathcal{O})} ds \\
&\quad + \epsilon C(C_0) \int_{x_0}^x H(y,s)\|(\omega^k, \nabla \omega^k)(y,s)\|^2_{H^{m-1}(\mathcal{O})} ds.
\end{aligned}$$

As we did before, we choose in the above inequality $\Lambda = \Lambda_1(h, c_0, \bar{c}_0, \epsilon)$ large enough so as to get rid of the terms containing $\|\partial_y \omega^{k+1}(y, s)\|^2_{H^{m-1}(\mathscr{O})}$. We then obtain

$$
\begin{aligned}
H(y,x)\|\partial_y\omega^{k+1}(y,x)\|^2_{H^{m-1}(\mathscr{O})} &\leq H(y,x_0)\|\partial_y\omega_0^{-,k+1}(y)\|^2_{H^{m-1}(\mathscr{O})} \\
&\quad + c_3(h,c_0,\bar{c}_0)\int_{x_0}^x H(y,s)\|\nabla_x\omega^{k+1}(y,s)\|^2_{H^m(\mathscr{O})}ds \\
&\quad + \epsilon C(C_0)\int_{x_0}^x H(y,s)\|(\omega^k,\nabla\omega^k)(y,s)\|^2_{H^{m-1}(\mathscr{O})}ds .
\end{aligned}
\tag{B.14}
$$

Then according to Lemma B.2, we estimate the terms containing $\|\nabla_x\omega^{k+1}(y,s)\|^2_{H^m(\mathscr{O})}$ and $\|(\omega^k, \nabla_x\omega^k)(y,s)\|^2_{H^{m-1}(\mathscr{O})}$ by using inequality (B.12) twice: first as it is written and secondly by replacing in that inequality $k$ with $k-1$ (which remains true according to Remark B.3):

$$
\begin{aligned}
H(y,x)\|\partial_y\omega^{k+1}(y,x)\|^2_{H^{m-1}(\mathscr{O})} &\leq H(y,x_0)\|\partial_y\omega_0^{-,k+1}(y)\|^2_{H^{m-1}(\mathscr{O})} + 2\epsilon C(C_0)\hat{C}(y,x) \\
&\quad + 2c_3(h,c_0,\bar{c}_0)\hat{C}(y,x) \\
&\quad + \epsilon C(C_0)\int_{x_0}^x H(y,s)\|\partial_y\omega^k(y,s)\|^2_{H^{m-1}(\mathscr{O})}ds .
\end{aligned}
\tag{B.15}
$$

We then integrate (B.15) with respect to $x$ on $[x_0, v]$ for any $v \in [x_0, 0[$

$$
\begin{aligned}
\int_{x_0}^v \overline{H}(x)\|\partial_y\omega^{k+1}(y,x)\|^2_{H^{m-1}(\mathscr{O})}dx &\leq |x_0|\left(\overline{H}(x_0)\|\partial_y\omega_0^{-,k+1}(y)\|^2_{H^{m-1}(\mathscr{O})} + 2e^{\lambda y}\epsilon C(C_0)\hat{C}(y,0)\right) \\
&\quad + 2c_2(h,c_0,\bar{c}_0)\hat{C}(y,0)|x_0|e^{\lambda y} \\
&\quad + \epsilon C(C_0)\int_{x_0}^v\int_{x_0}^x \overline{H}(s)\|\partial_x\omega^k(y,s)\|^2_{H^{m-1}(\mathscr{O})}dsdx ;
\end{aligned}
\tag{B.16}
$$

recall $\overline{H}(x) = |x|^{-2\alpha}e^{-\Lambda x}$. Let

$$
\begin{aligned}
\widetilde{C}(y_0, 0) &= \sup_{k\in\mathbb{N},\, y\in[0,y_0]}\left\{|x_0|\left(\overline{H}(x_0)\|\partial_y\omega_0^{-,k+1}(y)\|^2_{H^{m-1}(\mathscr{O})} + 2e^{\lambda y}\epsilon C(C_0)\hat{C}(y,0)\right) + 2c_2(h,c_0,\bar{c}_0)\hat{C}(y_0,0)|x_0|e^{\Lambda y}\right\} \\
&\quad + \frac{1}{2}\sup_{y\in[0,y_0]}\int_{x_0}^0 \overline{H}(x)\|\partial_y\omega^0(y,x)\|^2_{H^{m-1}(\mathscr{O})}dx .
\end{aligned}
\tag{B.17}
$$

Again we need to prove that $\widetilde{C}(y_0, 0)$ is a finite quantity. For all $\gamma \in \mathbb{N}^{n-1}$ such that $|\gamma| \leq m-1$, we have the following trivial identity

$$
\partial_y|\partial_\theta^\gamma\partial_y\omega_0^{-,k+1}|^2 = 2\partial_y\partial_\theta^\gamma\partial_y^2\omega_0^{-,k+1} \cdot \partial_\theta^\gamma\partial_y\omega_0^{-,k+1} \leq |\partial_\theta^\gamma\partial_y^2\omega_0^{-,k+1}|^2 + |\partial_\theta^\gamma\partial_y\omega_0^{-,k+1}|^2 .
$$

Integrating with respect to $y$ on the interval $[0, y_0]$ gives

$$
|\partial_\theta^\gamma\partial_y\omega_0^{-,k+1}(y,\theta)|^2 = |\partial_\theta^\gamma\partial_y\omega_0^{-,k+1}(0,\theta)|^2 + \int_0^{y_0}\left(|\partial_\theta^\gamma\partial_y^2\omega_0^{-,k+1}(s,\theta)|^2 + |\partial_\theta^\gamma\partial_y\omega_0^{-,k+1}(s,\theta)|^2\right)ds.
$$

We integrate this new identity now with respect to the angular variables on $\mathscr{O}$ and obtain that

$$
\|\partial_y\omega_0^{-,k+1}(y)\|^2_{H^{m-1}(\mathscr{O})} \leq \|\partial_y\omega_0^{-,k+1}(0)\|^2_{H^{m-1}(\mathscr{O})} + 2\|\omega_0^{-,k+1}\|^2_{H^{m+1}(\mathcal{C}^-)} .
$$

By the trace theorem (recall $\partial\mathcal{C}^- = (\{0\} \times \mathscr{O}) \cup (\{y_0\} \times \mathscr{O})$):

$$
\|\partial_y\omega_0^{-,k+1}(0)\|^2_{H^{m-1}(\mathscr{O})} \leq c\|\partial_y\omega_0^{-,k+1}\|^2_{H^m(\mathcal{C}^-)} \leq c\|\omega_0^{-,k+1}\|^2_{H^{m+1}(\mathcal{C}^-)} .
$$

It then follows that
$$\|\partial_y \omega_0^{-,k+1}(y)\|^2_{H^{m-1}(\mathscr{O})} \leq c \|\omega_0^{-,k+1}\|^2_{H^{m+1}(\mathcal{C}^-)} < \infty \,. \tag{B.18}$$

On the other hand, by the Equation 3.22 page 12 which defined $\omega_0$, we have
$$\begin{aligned}
\int_{x_0}^0 \overline{H}(x) \|\partial_y \omega^0(y,x)\|^2_{H^{m-1}(\mathscr{O})} dx &\leq \|\partial_y \omega_0^{-,0}(y)\|^2_{H^{m-1}(\mathscr{O})} \int_{x_0}^0 |x|^{-2\alpha} e^{-\Lambda x} dx \\
&\leq C \quad \text{independently of } y \,.
\end{aligned} \tag{B.19}$$

The estimates (B.18) and (B.19) prove that (B.17) defines a finite quantity. By the definition of $\widetilde{C}(y_0, 0)$, inequality (B.16) can be rewritten as:
$$\int_{x_0}^v \overline{H}(x) \|\partial_y \omega^{k+1}(y,x)\|^2_{H^{m-1}(\mathscr{O})} dx \leq \widetilde{C}(y_0, 0) + \epsilon C(C_0) \int_{x_0}^v \int_{x_0}^x \overline{H}(s) \|\partial_y \omega^k(y,s)\|^2_{H^{m-1}(\mathscr{O})} ds\, dx \,.$$

Suppose that $\forall v \in [x_0, 0[$,
$$\int_{x_0}^v \overline{H}(x) \|\partial_y \omega^k(y,x)\|^2_{H^{m-1}(\mathscr{O})} dx \leq 2\widetilde{C}(y_0, 0). \tag{B.20}$$

Then inequality (B) gives:
$$\begin{aligned}
\int_{x_0}^v \overline{H}(x) \|\partial_y \omega^{k+1}(y,x)\|^2_{H^{m-1}(\mathscr{O})} dx &\leq \widetilde{C}(y_0, 0) + 2\epsilon C(C_0) \widetilde{C}(u_0, 0) |x_0| \\
&\leq 2\widetilde{C}(y_0, 0), \quad \text{for } \epsilon \text{ small enough.}
\end{aligned}$$

Note that from the definition of the constant $\widetilde{C}(y_0, 0)$, inequality (B.20) remains true when $k = 0$ and then one can conclude that it holds for any integer $k \in \mathbb{N}$. Inequality (B.15) then implies:
$$|x|^{-2\alpha} \|\partial_y \omega^{k+1}(y,x)\|^2_{H^{m-1}(\mathscr{O})} \leq C_1(c_0, \bar{c}_0, h, C_0, \Lambda_0, \Lambda_1, \epsilon, \hat{c}, \widetilde{C}) \,. \tag{B.21}$$

In order to obtain the analog of (B.21) with instead $\partial_x \omega^{k+1}$, we repeat the previous argument. Once more we differentiate with $\partial_\theta^\gamma$ the equation satisfies by $\omega^{k+1}$ and multiply the resulting equation by $H \partial_\theta^\gamma \partial_x \omega^{k+1}$ and obtain
$$\begin{aligned}
\partial_y \left( H (\partial_x \partial_\theta^\gamma \omega^{k+1})^2 \right) &= (\partial_y H + H \frac{n-1}{2\rho})(\partial_x \partial_\theta^\gamma \omega^{k+1})^2 - H \frac{n-1}{2\rho} \partial_y \partial_\theta^\gamma \omega^{k+1} \partial_x \partial_\theta^\gamma \omega^{k+1} \\
&\quad + \partial_x \partial_\theta^\gamma \omega^{k+1} \sum_{\gamma_1 + \gamma_2 = \gamma} H \frac{\partial_\theta^{\gamma_1} h^{AB}}{2\rho^2} \partial_\theta^{\gamma_2} \partial_A \partial_B \omega^{k+1} \\
&\quad - \partial_x \partial_\theta^\gamma \omega^{k+1} \sum_{\gamma_1 + \gamma_2 = \gamma} H \frac{\partial_\theta^{\gamma_1} \Gamma^B}{2\rho^2} \partial_\theta^{\gamma_2} \partial_B \omega^{k+1} \\
&\quad - \frac{1}{2} |x|^{-\frac{n+3}{2}} H \partial_x \partial_\theta^\gamma \omega^{k+1} \partial_\theta^\gamma G^k(\ldots) \,.
\end{aligned}$$

Then, we integrate on $[0, y] \times \{x\} \times \mathscr{O}$, and obtain for $n \geq 1 + \frac{4}{r-1} - 2\alpha$ via Stokes theorem
$$\begin{aligned}
H(y,x) \|\partial_x \omega^{k+1}(y,x)\|^2_{H^{m-1}(\mathscr{O})} &\leq H(0,x) \|\partial_x \omega_0^{+,k+1}(x)\|^2_{H^{m-1}(\mathscr{O})} \\
&\quad + (c_4(h, c_0, \bar{c}_0) + c(\epsilon) - \Lambda) \int_0^y H(s,x) \|\partial_x \omega^{k+1}(s,x)\|^2_{H^{m-1}(\mathscr{O})} ds \\
&\quad + c_5(h, c_0, \bar{c}_0) \int_0^y H(s,x) \|\nabla_y \omega^{k+1}(s,x)\|^2_{H^m(\mathscr{O})} ds \\
&\quad + \epsilon C(C_0) \int_0^y H(s,x) \|(\omega^k, \nabla \omega^k)(s,x)\|^2_{H^{m-1}(\mathscr{O})} ds
\end{aligned}$$

As we did previously, we choose in this inequality $\Lambda = \Lambda_2(h, c_0, \bar{c}_0, \epsilon)$ large enough so as to get rid of the terms with $\|\partial_x \omega^{k+1}(y,s)\|^2_{H^{m-1}(\mathscr{O})}$ and we obtain:

$$
\begin{aligned}
H(y,x)\|\partial_x\omega^{k+1}(y,x)\|^2_{H^{m-1}(\mathscr{O})} &\leq H(0,x)\|\partial_x\omega_0^{+,k+1}(x)\|^2_{H^{m-1}(\mathscr{O})} \\
&\quad + c_5(h,c_0,\bar{c}_0)\int_0^y H(s,x)\|\nabla_y\omega^{k+1}(s,x)\|^2_{H^m(\mathscr{O})}ds \\
&\quad + \epsilon C(C_0)\int_0^y H(s,x)\|(\omega^k, \nabla\omega^k)(s,x)\|^2_{H^{m-1}(\mathscr{O})}ds .
\end{aligned}
\tag{B.22}
$$

We estimate the quantities $\|\nabla_y\omega^{k+1}(y,s)\|^2_{H^m}$ and $\|(\omega^k, \nabla_y\omega^k)(y,s)\|^2_{H^{m-1}(\mathscr{O})}$ using inequality (B.12):

$$
\begin{aligned}
H(y,x)\|\partial_x\omega^{k+1}(y,x)\|^2_{H^{m-1}(\mathscr{O})} &\leq H(0,x)\|\partial_x\omega_0^{+,k+1}(x)\|^2_{H^{m-1}(\mathscr{O})} + 2\epsilon C(C_0)\hat{C}(y_0,0) \\
&\quad + 2c_5(h,c_0,\bar{c}_0)\hat{C}(y_0,0) \\
&\quad + \epsilon C(C_0)\int_0^y H(s,x)\|\partial_x\omega^k(s,x)\|^2_{H^{m-1}(\mathscr{O})}ds .
\end{aligned}
\tag{B.23}
$$

We integrate in $y$, and obtain that for any $u \in [0, y_0]$,

$$
\begin{aligned}
\int_0^u \widetilde{H}(y,x)\|\partial_x\omega^{k+1}(y,x)\|^2_{H^{m-1}(\mathscr{O})}dy &\leq y_0\widetilde{H}(0,x)\|\partial_x\omega_0^{+,k+1}(x)\|^2_{H^{m-1}(\mathscr{O})} \\
&\quad + 2y_0 e^{\Lambda x}\epsilon C(C_0)\hat{C}(y_0,0) + 2c_5(h,c_0,\bar{c}_0)\hat{C}(y_0,0)y_0 e^{\Lambda x} \\
&\quad + \epsilon C(C_0)\int_0^u\int_0^y \widetilde{H}(s,x)\|\partial_x\omega^k(y,s)\|^2_{H^{m-1}(\mathscr{O})}dsdy ,
\end{aligned}
\tag{B.24}
$$

where $\widetilde{H}(y,x) = |x|^{-2\alpha}e^{-\Lambda y}$. Now, we define a new constant $\check{C}(y_0, 0)$ as

$$
\begin{aligned}
\check{C}(y_0, 0) &= \sup_{k\in\mathbb{N},\, x\in[x_0, 0[} \left\{ y_0\left(\widetilde{H}(0,x)\|\partial_x\omega_0^{+,k+1}(x)\|^2_{H^{m-1}(\mathscr{O})} + 2e^{\Lambda x}\epsilon C(C_0)\hat{C}(y_0,0)\right) \right. \\
&\qquad\qquad\qquad \left. + 2c_5(h,c_0,\bar{c}_0)\hat{C}(y_0,0)|x_0|e^{\Lambda x} \right\} \\
&\quad + \frac{1}{2}\sup_{x\in[x_0,0[}\int_0^{y_0}\widetilde{H}(y,x)\|\partial_x\omega^0(s,x)\|_{H^{m-1}(\mathscr{O})}.
\end{aligned}
\tag{B.25}
$$

As before let us prove that (B.25) is finite. By the Lemma 3.8 page 12

$$
\sup_{k\in\mathbb{N},\, x\in[x_0,0[} \widetilde{H}(0,x)\|\partial_x\omega_0^{+,k+1}(x)\|^2_{H^{m-1}(\mathscr{O})} < \infty .
\tag{B.26}
$$

Next, by the definition of $\omega_0$ given by (3.22) page 12, we have

$$
\begin{aligned}
\int_0^{y_0} \widetilde{H}(s,x)\|\partial_x\omega^0(s,x)\|^2_{H^{m-1}(\mathscr{O})}ds &\leq \|\partial_x\omega_0^{+,0}(x)\|^2_{H^{m-1}(\mathscr{O})}|x|^{-2\alpha}\int_0^{y_0}e^{-\Lambda s}ds \\
&\leq C \text{ independently of } x .
\end{aligned}
\tag{B.27}
$$

From the estimates (B.26) and (B.27) we obtain that (B.25) define a finite quantity. We thus obtain the following form of (B.24):

$$
\int_0^u \widetilde{H}(y,x)\|\partial_x\omega^{k+1}(y,x)\|^2_{H^{m-1}(\mathscr{O})}dy \leq \check{C} + \epsilon C(C_0)\int_0^u\int_0^y \widetilde{H}(s,x)\|\partial_x\omega^k(s,x)\|^2_{H^{m-1}(\mathscr{O})}dsdy .
\tag{B.28}
$$

Suppose again that $\forall u \in [0, y_0]$,

$$\int_0^u \widetilde{H}(y,x) \|\partial_x \omega^k(y,x)\|_{H^{m-1}(\mathscr{O})}^2 dy \leq 2\check{C}(y_0, 0). \tag{B.29}$$

Then inequality (B.28) gives:

$$\int_0^u \widetilde{H}(y,x) \|\partial_x \omega^{k+1}(y,x)\|_{H^{m-1}(\mathscr{O})}^2 dy \leq \check{C}(y_0, 0) + 2\epsilon C(C_0)\check{C}(y_0, 0) y_0$$
$$\leq 2\check{C}(y_0, 0), \quad \text{if } \epsilon \text{ small enough.}$$

By the definition of the constant $\check{C}(y_0, 0)$, (B.29) is satisfied for $k = 0$ and so it is for any integer $k \in \mathbb{N}$. Inequality (B.23) implies:

$$|x|^{-2\alpha} \|\partial_x \omega^{k+1}(y,x)\|_{H^{m-1}(\mathscr{O})}^2 \leq C_2(c_0, \bar{c}_0, h, C_0, \Lambda_0, \Lambda_2, \epsilon, \hat{C}, \check{C}). \tag{B.30}$$

It remains to control the $H^m(\mathscr{O})$ norms of $\omega^{k+1}$, that is to control its angular derivatives. Let $(y, x) \in [0, y_0] \times [x_0, 0[$, $\gamma \in \mathbb{N}^{n-1}$ such that $|\gamma| \leq m$

$$|x|^{-\alpha} e^{-\frac{\Lambda}{2} y} |\partial_\theta^\gamma \omega^{k+1}(y,x)| \leq |x|^{-\alpha} |\partial_\theta^\gamma \omega_0^{+,k+1}(x)| + \int_0^y |x|^{-\alpha} e^{-\frac{\Lambda}{2} s} \partial_y \partial_\theta^\gamma \omega^{k+1}(s,x)| ds.$$

It then follows that:

$$|x|^{-2\alpha} |e^{-\Lambda y} \partial_\theta^\gamma \omega^{k+1}(y,x)|^2 \leq 2\left(|x|^{-2\alpha} \partial_\theta^\gamma \omega_0^{+,k+1}(x)|^2 + y_0 \int_0^y |x|^{-2\alpha} |e^{-\Lambda s} \partial_\theta^\gamma \partial_y \omega^{k+1}(s,x)|^2 ds\right).$$

By integration we then obtain:

$$|x|^{-2\alpha} e^{-\Lambda y} \|\omega^{k+1}(y,x)\|_{H^m(\mathscr{O})}^2 \leq 2\left(|x|^{-2\alpha} \|\omega_0^{+,k+1}(x)\|_{H^m(\mathscr{O})}^2 + y_0 \int_0^y |x|^{-2\alpha} e^{-\Lambda s} \|\partial_y \omega^{k+1}(s,x)\|_{H^m(\mathscr{O})}^2 ds\right)$$
$$\leq C_3(\hat{C}) \quad \text{see B.12}.$$

We have proved the following Lemma:

**Lemma B.4** *Let $m \in \mathbb{N}^*$. If $n \geq 1 + \frac{4}{r-1} - 2\alpha$, then there exists a positive constant $C_4 = C_4(c_0, \bar{c}_0, h, \Lambda_0, \Lambda_1, \Lambda_2, \epsilon_0)$ such that:*

$$\sup_{(y,x) \in [0, y_0] \times [x_0, 0[} x^{-\alpha} \|(\omega^{k+1}, \nabla \omega^{k+1})\|_{H^{m-1}(\mathscr{O})} \leq C_4. \tag{B.31}$$

Now to prove that (B.1) holds for $k = j + 1$ we are going to show that it suffices to replace in (B.31) $y_0$ with a certain $u_*$ sufficiently small. Let $j_0 \in \mathbb{N}^*$. if

$$m - 1 > \frac{n-1}{2} + j_0,$$

then from (B.31) and the Sobolev embedding theorem, for all $(y,x) \in [0, y_0] \times [x_0, 0[$, we have:

$$|x|^{-\alpha} \|(\omega^{k+1}, \nabla \omega^{k+1})(y,x)\|_{C^{j_0}(\mathscr{O})} \leq C. \tag{B.32}$$

It then follows from the differential equation satisfies by $\omega^{k+1}$ that

$$\|\partial_y(|x|^{-\alpha} \partial_x \omega^{k+1})(y,x)\|_{C^{j_0-1}(\mathscr{O})} \leq C, \tag{B.33}$$
$$\|x\partial_x(|x|^{-\alpha} \partial_y \omega^{k+1})(y,x)\|_{C^{j_0-1}(\mathscr{O})} \leq C. \tag{B.34}$$

Integrating (B.33) in $y$ from $(0,x)$ to $(y,x)$ we find that for $j_0 \geq 2$,

$$\begin{aligned} |x|^{-\alpha}\|\partial_x\omega^{k+1}(y,x)\|_{C^1(\mathscr{O})} &\leq |x|^{-\alpha}\|\partial_x\omega_0^{+,k+1}(0,x)\|_{C^1(\mathscr{O})} + Cy \\ &\leq C_0 + Cy \\ &\leq 2C_0 \quad \text{if} \quad y \leq u_1 \ . \end{aligned} \qquad (\text{B.35})$$

Note that inequality (B.35) shall be read as a first condition in the determination of $u_*$. Further, to control $|x|^{-\alpha}\|\partial_y\omega^{k+1}(y,x)\|_{C^1(\mathscr{O})}$, we $y$–differentiate the differential equation satisfied by $\omega^{k+1}$. We write here

$$G(z, x^{\frac{n-1}{2}}(\omega^k, \partial_x\omega^k, \nabla_y\omega^k) = G(z, x^{\frac{n-1}{2}}(p_1, p_2, p_3))$$

so that $\partial_y G^k(\ldots)$ reads

$$\begin{aligned} \partial_y G^k(\ldots) &= (\partial_y G)^k(\ldots) + |x|^{\frac{n-1}{2}} \partial_y\omega^k \left(\frac{\partial G}{\partial p_1}\right)^k (\ldots) \\ &\quad + |x|^{\frac{n-1}{2}} \partial_y\partial_x\omega^k \left(\frac{\partial G}{\partial p_2}\right)^k (\ldots) + |x|^{\frac{n-1}{2}} \partial_y\nabla_y\omega^k \left(\frac{\partial G}{\partial p_3}\right)^k (\ldots) \ . \end{aligned}$$

The differentiated equation reads

$$\partial_x(\partial_y\nabla_y\omega^{k+1}) = \xi\partial_y\nabla_y\omega^{k+1} + \psi^k \partial_y\nabla_y\omega^k + \varphi^k \ , \qquad (\text{B.36})$$

where we have set

$$4\xi = \left(-\frac{n-1}{\rho} \ , \ 0\right), \qquad 4\psi^k = -|x|^{-2}\left(\frac{\partial G}{\partial p_3}\right)^k(\ldots) \ ,$$

and where the components of $\varphi^k$ are given by

$$\begin{aligned} 4\varphi_y^k(y,x) &= \frac{n-1}{\rho}\partial_x\partial_y\omega^{k+1}(y,x) + \frac{h^{AB}}{\rho^3}\partial_A\partial_B\omega^{k+1}(y,x) + \frac{h^{AB}}{\rho^2}\partial_y\partial_A\partial_B\omega^{k+1}(y,x) \\ &\quad + \frac{n-1}{2\rho^2}(\partial_x - \partial_y)\omega^{k+1} - \frac{\Gamma^B}{\rho^2}\partial_y\partial_B\omega^{k+1} - \frac{\Gamma^B}{\rho^3}\partial_B\omega^{k+1} - (-x)^{\frac{n+3}{2}}(\partial_y G)^k(\ldots) \\ &\quad - x^{-2}\left(\partial_y\omega^k\left(\frac{\partial G}{\partial p_1}\right)^k(\ldots) + \partial_y\partial_x\omega^k\left(\frac{\partial G}{\partial p_2}\right)^k(\ldots)\right) \ , \\ \varphi_A^k(y,x) &= \partial_y\partial_x\partial_A\omega^{k+1} \ . \end{aligned}$$

By hypothesis (3.9)) and the induction assumption (B.1), we have the following estimate on $\psi^k(y,x)$ for all $(y,x) \in [0,y_0] \times [x_0, 0[$:

$$\begin{aligned} \|\psi^k(y,x)\|_{C^{j_0}} &\leq C(|x|^{-\alpha}\|(\omega^k, \nabla\omega^k)(y,x)\|_{L^\infty(\mathscr{O})})|x|^{-2+(r-1)(\frac{n-1}{2}+\alpha)} \\ &\leq C \quad \text{if} \quad n \geq 1 + \frac{4}{r-1} - 2\alpha \ . \end{aligned} \qquad (\text{B.37})$$

By using simultaneously (B.32) and (B.33), for all $(y,x) \in [0,y_0] \times [x_0, 0[$ we have:

$$|x|^{-\alpha}\|\varphi^k(y,x)\|_{C^{j_0-2}(\mathscr{O})} < C \quad \text{if} \quad n \geq 1 + \frac{4}{r-1} - 2\alpha \ . \qquad (\text{B.38})$$

On the other hand we have:

$$\begin{aligned} \partial_x\left(\overline{H}(x)|\partial_y\nabla_y\omega^{k+1}|^2\right) &= (2\alpha|x|^{-1} - \Lambda)\overline{H}(x)|\partial_y\nabla_y\omega^{k+1}|^2 + 2\overline{H}(x)\partial_x\partial_y\nabla_y\omega^{k+1}.\partial_y\nabla_y\omega^{k+1} \\ &\leq -\Lambda\overline{H}(x)|\partial_y\nabla_y\omega^{k+1}|^2 + 2\overline{H}(x)\partial_x\partial_y\nabla_y\omega^{k+1}.\partial_y\nabla_y\omega^{k+1} \ . \end{aligned}$$

Considering now inequalities (B.37) and (B.38), we deduce that for $j_0 \geq 2$,

$$\begin{aligned}\partial_x \left(\overline{H}(x)|\partial_y \nabla_y \omega^{k+1}|^2\right) &\leq -\Lambda \overline{H}(x)|\partial_y \nabla_y \omega^{k+1}|^2 + 2\overline{H}(x)\partial_y \nabla_y \omega^{k+1} \left(\xi \partial_y \nabla_y \omega^{k+1} + \psi^k \partial_y \nabla_y \omega^k + \varphi^k\right) \\ &\leq \overline{H}(x)\left((c + \frac{C}{\epsilon} - \Lambda)|\partial_y \nabla_y \omega^{k+1}|^2 + \epsilon C|\partial_y \nabla_y \omega^k|^2\right) + C|x|^{-\alpha} e^{-\Lambda x}.\end{aligned}$$

Choosing $\Lambda$ large enough, we have thus proved the following inequality

$$\partial_x \left(\overline{H}(x)|\partial_y \nabla_y \omega^{k+1}|^2\right) \leq \epsilon \overline{H}(x) C |\partial_y \nabla_y \omega^k|^2 + C|x|^{-\alpha} e^{-\Lambda x}, \tag{B.39}$$

which is then integrated in $x$ to obtain:

$$\overline{H}(x)|\partial_y \nabla_y \omega^{k+1}|^2(y,x) \leq \overline{H}(x_0)|\partial_y \nabla_y \omega^{k+1}|^2(y,x_0) + C \int_{x_0}^x e^{-\Lambda s} \left(\epsilon |s|^{-2\alpha}|\partial_y \nabla_y \omega^k|^2 + |s|^{-\alpha}\right) ds.$$

Equivalently this reads

$$\overline{H}(x)|\partial_y \nabla_y \omega^{k+1}|^2(y,x) \leq \overline{H}(x_0)|\partial_y \nabla_y \omega_0^{-,k+1}|^2(y,x_0) + C \int_{x_0}^x e^{-\Lambda s} \left(\epsilon |s|^{-2\alpha}|\partial_y \nabla_y \omega^k|^2 + |s|^{-\alpha}\right) ds$$

As we did many times before, for a convenient choice of $\epsilon$, we obtain the estimate:

$$\begin{aligned}|x|^{-2\alpha}|\partial_y \nabla_y \omega^{k+1}|^2(y,x) &\leq 2\left(C \int_{x_0}^0 |s|^{-\alpha} e^{-\Lambda s} ds + \sup_{k \in \mathbb{N},\ y \in [0,y_0]} \overline{H}(x_0)|\partial_y \nabla_y \omega_0^{-,k}|^2(y)\right) \\ &\leq c(x_0,\alpha)\left(1 + \sup_{k \in \mathbb{N}} \|\omega_0^{-,k+1}\|^2_{C^2([0,y_0] \times \mathscr{O})}\right) \\ &\leq c(x_0,\alpha)\left(1 + \sup_{k \in \mathbb{N}} \|\omega_0^{-,k+1}\|^2_{H^{m+1}(\mathcal{C}^-)}\right) < \infty\end{aligned}$$

i.e.

$$|x|^{-\alpha}|\partial_y \nabla_y \omega^{k+1}|(y,x) < C. \tag{B.40}$$

The same procedure can exactly be repeated by using instead the angular derivatives of $\omega^{k+1}$ leading to

$$|x|^{-\alpha}\|\partial_y \nabla_y \omega^{k+1}(y,x)\|_{C^{j_0-2}} < C. \tag{B.41}$$

**Remark B.5** The bound (B.41) follows from a control on the coefficients $\xi$, $\psi^k$ and $\varphi^k$ in (B.36) after have y-differentiated the equation satisfied by $\omega^{k+1}$. We notice that, if instead we have $x\partial_x-$ differentiated the same equation we should have been led to the bound:

$$|x|^{-\alpha}\|(x\partial_x)\nabla_x \omega^{k+1})(y,x)\|_{C^{j_0-2}} < C. \tag{B.42}$$

Now we integrate once more in $y$ Equation B.41 and obtain that for $j_0 \geq 3$:

$$\begin{aligned}|x|^{-\alpha}\|\nabla_y \omega^{k+1}(y,x)\|_{C^1(\mathscr{O})} &\leq |x|^{-\alpha}\|\nabla_y \omega^{k+1}(0,x)\|_{C^1} + Cy \\ &\leq 2C_0 \quad \text{pour} \quad y < u_2.\end{aligned} \tag{B.43}$$

In order to finish the proof of Lemma 3.9, we write:

$$\begin{aligned}|x|^{-\alpha}\|\omega^{k+1}(y,x)\|_{C^1(\mathscr{O})} &\leq |x|^{-\alpha}\|\omega_0^{+,k+1}(x)\|_{C^1(\mathscr{O})} + \int_0^y |x|^{-\alpha}\|\partial_y \omega^{k+1}(s,x)\|_{C^1(\mathscr{O})} ds \\ &\leq C_0 + 2C_0 y \\ &\leq 2C_0 \quad \text{pour} \quad y < u_3.\end{aligned} \tag{B.44}$$

We define $u_*$ as

$$u_* = \min\{u_0, u_1, u_2, u_3\}$$

and inequalities (B.35), (B.43), (B.44) allow us to write:

$$\sup_{(y,x) \in [0,u_*] \times [x_0,0[} |x|^{-\alpha}\|(\omega^{k+1}, \nabla \omega^{k+1})(y,x)\|_{W^{1,\infty}} < 2C_0.$$

This completes the proof of Lemma 3.9. $\square$

# C Proof of Lemma 3.12 page 14

First we will prove (3.28) and secondly, we will show that (3.27) is actually a consequence of (3.28). We proceed by induction on $k$. Set

$$
\begin{aligned}
\bar{C}_0 &:= \sup_{(y,x)\in[0,u_*]\times[x_0,0[} |x|^{-2\alpha}\|(\partial_x\omega^0,\partial_y\omega^0)(y,x)\|^2_{H^{m-1}(\mathscr{O})} < \infty\ , \\
\bar{C}_1 &:= \sup_{k\in\mathbb{N},\ (y,x)\in[0,u_*]\times[x_0,0[} \Big\{ H(y,x_0)\|\partial_y\omega_0^{-,k+1}(y)\|^2_{H^{m-1}(\mathscr{O})} \\
&\qquad\qquad\qquad\qquad\qquad +H(0,x)\|\partial_x\omega_0^{+,k+1}(x)\|^2_{H^{m-1}(\mathscr{O})} \Big\} \\
&\quad +2\epsilon C(C_0)\hat{C}(y_0,0) + 2c_3(h,c_0,\bar{c}_0)\hat{C}(y_0,0) + 2\epsilon C(C_0)\hat{C}(y_0,0) \\
&\quad +2c_5(h,c_0,\bar{c}_0)\hat{C}(y_0,0) < \infty\ ,
\end{aligned} \tag{C.1}
$$

(the constant $\hat{C}(y_0,0)$ is defined in Equation B.8 page 26) and suppose that

$$\sup_{(y,x)\in[0,u_*]\times[x_0,0[} H(y,x)\|(\partial_x\omega^k,\partial_y\omega^k)(y,x)\|^2_{H^{m-1}(\mathscr{O})} < 2(\bar{C}_0+\bar{C}_1)\ .$$

Let us show that this remains true when we replace $k$ with $k+1$. Adding inequalities (B.15) and (B.23) leads to:

$$
\begin{aligned}
H(y,x)\|(\partial_y\omega^{k+1},\partial_x\omega^{k+1})(y,x)\|^2_{H^{m-1}(\mathscr{O})} &\leq \bar{C}_1 + \epsilon C(C_0)\int_{x_0}^x H(y,s)\|\partial_y\omega^k(y,s)\|^2_{H^{m-1}(\mathscr{O})}ds \\
&\quad +\epsilon C(C_0)\int_0^y H(s,x)\|\partial_x\omega^k(s,x)\|^2_{H^{m-1}(\mathscr{O})}ds \\
&\leq \bar{C}_1 + 2\epsilon C(C_0)(\bar{C}_0+\bar{C}_1)(y_0+|x_0|) \\
&\leq 2(\bar{C}_0+\bar{C}_1)\ ,\quad \text{even if it means redefining } \epsilon\ . \quad (C.2)
\end{aligned}
$$

Since $\forall(y,x)\in[0,y_0]\times[x_0,0[,\ e^{\Lambda(x+y)}\leq e^{\Lambda y_0}$, it then suffices to set

$$M_2 := (2e^{\Lambda y_0}(\bar{C}_0+\bar{C}_1))^{1/2}\ .$$

In order to obtain the uniform control (3.27), we repeat the argument leading to the proof of Lemma 3.11. For all $(y,x)\in[0,u_*]\times[x_0,0[$, we have:

$$
\begin{aligned}
\|\omega^k(y,x)\|_{H^{m-1}(\mathscr{O})} &\leq \|\omega_0^{-,k}(y)\|_{H^{m-1}(\mathscr{O})} + \int_{x_0}^x \|\partial_x\omega^k(y,s)\|_{H^{m-1}(\mathscr{O})}ds \\
&\leq \|\omega_0^{-,k}(y)\|_{H^{m-1}(\mathscr{O})} + \int_{x_0}^x |s|^\alpha|s|^{-\alpha}\|\partial_x\omega^k(s,x)\|_{H^{m-1}(\mathscr{O})}ds \\
&\leq \underbrace{\sup_{k\in\mathbb{N},\ y\in[0,u_*]}\|\omega_0^{-,k}(y)\|_{H^{m-1}(\mathscr{O})} + M_2\int_{x_0}^0 |s|^\alpha ds}_{:=M_1}\ . \quad (C.3)
\end{aligned}
$$

$\square$

# D  Proof of Lemma 3.13, page 14

We apply Proposition 3.6 with $\omega = \delta\omega^k$, $u \in [0, u_*]$, and $v \in [x_0, 0[$. We have:

$$e^{-\Lambda u}\|(\delta\omega^k, \delta(\nabla_x\omega^k))\overline{H}^{\frac{1}{2}}\|^2_{L^2([x_0,v[\times\mathscr{O})} + e^{-\Lambda v}\|(\delta\omega^k, \delta(\nabla_y\omega^k))\widetilde{H}^{\frac{1}{2}}\|^2_{L^2([0,u]\times\mathscr{O})} \leq$$
$$\|(\delta\omega^k, \delta(\nabla_x\omega^k))\overline{H}^{\frac{1}{2}}\|^2_{L^2([x_0,v[\times\mathscr{O})} + e^{-\Lambda x_0}\|(\delta\omega^k, \delta(\nabla_y\omega^k))\widetilde{H}^{\frac{1}{2}}\|^2_{L^2([0,u]\times\mathscr{O})}$$
$$+(c_1(c_0, \bar{c}_0, n, h) - 2\Lambda)\int_0^u \int_{x_0}^v H(y,x)\|(\delta\omega^k, \delta\nabla\omega^k)\|^2_{L^2(\mathscr{O})}dxdy$$
$$+\frac{1}{c_0}\int_{D_{u,v}} |L^\ell[\delta\omega^k]|\, dsdyd\nu\, . \tag{D.1}$$

Recall that:
$$L^\ell[\delta\omega^k] = |x|^{-\frac{n+3}{2}} H(\delta\partial_x\omega^k + \delta\partial_y\omega^k)\left(G^k(\ldots) - G^{k-1}(\ldots)\right)\, .$$

We have:
$$|x|^{-\frac{n+3}{2}} H\left(G^k(\ldots) - G^{k-1}(\ldots)\right) = H\int_0^1 |x|^{-2}\xi^k(t,y,x)dt$$

with

$$\xi^k(t,y,x) = \delta\omega^{k-1}\frac{\partial G}{\partial p}\left(z, t|x|^{\frac{n-1}{2}}(\omega^k, \nabla\omega^k) + (1-t)|x|^{\frac{n-1}{2}}(\omega^{k-1}, \nabla\omega^{k-1})\right)$$
$$+\delta\nabla\omega^{k-1}\frac{\partial G}{\partial q}\left(z, t|x|^{\frac{n-1}{2}}(\omega^k, \nabla\omega^k) + (1-t)|x|^{\frac{n-1}{2}}(\omega^{k-1}, \nabla\omega^{k-1})\right)\, .$$

Using once more hypothesis (3.9), one is led to the following estimate

$$|x|^{-2}|\xi(t,y,x)| \leq C_5\, |x|^{-2+(r-1)(\alpha+\frac{n-1}{2})}\left(|\delta\omega^{k-1}| + |\delta\nabla\omega^{k-1}|\right)$$
$$\leq C_5\left(|\delta\omega^{k-1}| + |\delta\nabla\omega^{k-1}|\right);\quad \text{if}\quad n \geq 1 + \frac{4}{r-1} - 2\alpha\, .$$

We should point out that the constant $C_5$ depends on the quantity

$$\sup_{k\in\mathbb{N},\, (y,x)\in[0,u_*]\times[x_0,0[} |x|^{-\alpha}\|(\omega^k, \nabla\omega^k)\|_{L^\infty(\mathscr{O})}\, ,$$

which does neither depend upon $\Lambda$ nor on $k$. We then obtain that if $n \geq 1 + \frac{4}{r-1} - 2\alpha$,

$$\int_{D_{u,v}} |L^\ell[\delta\omega^k]| \leq C_5 \int_0^u \int_{x_0}^v \left(\|\partial_x\delta\omega^k(y,x)\|^2_{L^2(\mathscr{O})} + \|\partial_y\delta\omega^k\|^2_{L^2(\mathscr{O})} + \|(\delta\omega^{k-1}, \nabla\delta\omega^{k-1})\|^2_{L^2(\mathscr{O})}\right)H(y,x)dxdy\, ;$$

and inequality (D.1) implies

$$e^{-\Lambda u}\|(\delta\omega^k, \delta(\nabla_x\omega^k))\overline{H}^{\frac{1}{2}}\|^2_{L^2([x_0,v[\times\mathscr{O})} + e^{-\Lambda v}\|(\delta\omega^k, \delta(\nabla_y\omega^k))\widetilde{H}^{\frac{1}{2}}\|^2_{L^2([0,u]\times\mathscr{O})} \leq$$
$$\|(\delta\omega^k, \delta(\nabla_x\omega^k))\overline{H}^{\frac{1}{2}}\|^2_{L^2([x_0,v[\times\mathscr{O})} + e^{-\Lambda x_0}\|(\delta\omega^k, \delta(\nabla_y\omega^k))\widetilde{H}^{\frac{1}{2}}\|^2_{L^2([0,u]\times\mathscr{O})}$$
$$+(c_1(c_0, \bar{c}_0, n, h) + C_5 - 2\Lambda)\int_0^u \int_{x_0}^v H(y,x)\|(\delta\omega^k, \delta\nabla\omega^k)\|^2_{L^2(\mathscr{O})}dxdy$$
$$+C_5 \int_0^u \int_{x_0}^v H(y,x)\|(\delta\omega^{k-1}, \nabla\delta\omega^{k-1})\|^2_{L^2(\mathscr{O})}dxdy\, . \tag{D.2}$$

From this inequality and by a convenient choice of $\Lambda$ we obtain:

$$e^{-\Lambda u}\|(\delta\omega^k,\delta(\nabla_x\omega^k))\overline{H}^{\frac{1}{2}}\|^2_{L^2([x_0,0[\times\mathscr{O})} + e^{-\Lambda v}\|(\delta\omega^k,\delta(\nabla_y\omega^k))\widetilde{H}^{\frac{1}{2}}\|^2_{L^2([0,u]\times\mathscr{O})} \leq$$
$$\|(\delta\omega_0^{+,k},\delta(\nabla_x\omega_0^{+,k}))\overline{H}^{\frac{1}{2}}\|^2_{L^2([x_0,0[\times\mathscr{O})} + e^{-\Lambda x_0}\|(\delta\omega_0^{-,k},\delta(\nabla_y\omega_0^{-,k}))\widetilde{H}^{\frac{1}{2}}\|^2_{L^2([0,u]\times\mathscr{O})}$$
$$+C_5\int_0^u\int_{x_0}^v H(y,x)\|(\delta\omega^{k-1},\nabla\delta\omega^{k-1})\|^2_{L^2(\mathscr{O})}dxdy\ . \tag{D.3}$$

We have:

$$\|(\delta\omega_0^{+,k},\delta(\nabla_x\omega_0^{+,k}))\overline{H}^{\frac{1}{2}}\|^2_{L^2([x_0,0[\times\mathscr{O})} + e^{-\Lambda x_0}\|(\delta\omega_0^{-,k},\delta(\nabla_y\omega_0^{-,k}))\widetilde{H}^{\frac{1}{2}}\|^2_{L^2([0,u]\times\mathscr{O})} =$$
$$\|e^{-\frac{\Lambda x}{2}}|x|^{-\alpha}(\delta\omega_0^{+,k},\delta(\nabla_x\omega_0^{+,k}))\|^2_{L^2([x_0,0[\times\mathscr{O})}$$
$$+|x_0|^{-\alpha}e^{-\Lambda x_0}\|e^{-\frac{1}{2}\Lambda y}(\delta\omega_0^{-,k},\delta(\nabla_y\omega_0^{-,k}))\|^2_{L^2([0,u]\times\mathscr{O})}$$
$$\leq\ c(\Lambda,x_0)\left(\||x|^{-\alpha}(\delta\omega_0^{+,k},\delta(\nabla_x\omega_0^{+,k}))\|^2_{L^2([x_0,0[\times\mathscr{O})} + \|(\delta\omega_0^{-,k},\delta(\nabla_y\omega_0^{-,k}))\|^2_{L^2([0,y_0]\times\mathscr{O})}\right)\ .$$

Since the sequences $\left(|x|^{-\alpha}(\omega_0^{+,k},\nabla_x\omega_0^{+,k})\right)_{k\in\mathbb{N}}$ and $\left(\omega_0^{-,k},\nabla_y\omega_0^{-,k}\right)_{k\in\mathbb{N}}$ are convergent respectively in the spaces $L^2([x_0,0[\times\mathscr{O})$ and $L^2([0,y_0]\times\mathscr{O})$, we know that

$$\lim_{k\to\infty}\left(\||x|^{-\alpha}(\delta\omega_0^{+,k},\delta(\nabla_x\omega_0^{+,k}))\|^2_{L^2([x_0,0[\times\mathscr{O})} + \|(\delta\omega_0^{-,k},\delta(\nabla_y\omega_0^{-,k}))\|^2_{L^2([0,y_0]\times\mathscr{O})}\right) = 0\ .$$

Therefore, $\forall i \in \mathbb{N}$, $\exists k_i \in \mathbb{N}$, such that

$$c(\Lambda,x_0)\left(\||x|^{-\alpha}(\delta\omega_0^{+,k_i},\delta(\nabla_x\omega_0^{+,k_i}))\|^2_{L^2([x_0,0[\times\mathscr{O})} + \|(\delta\omega_0^{-,k_i},\delta(\nabla_y\omega_0^{-,k_i}))\|^2_{L^2([0,y_0]\times\mathscr{O})}\right) \leq \frac{1}{2^i}. \tag{D.4}$$

We then write inequality (D.3) with instead the subsequence $(\omega^{k_i})_{i\in\mathbb{N}}$ which will be denoted again $(\omega^k)_{k\in\mathbb{N}}$ and one obtains:

$$e^{-\Lambda u}\|(\delta\omega^k,\delta(\nabla_x\omega^k))\overline{H}^{\frac{1}{2}}\|^2_{L^2([x_0,v[\times\mathscr{O})} + e^{-\Lambda v}\|(\delta\omega^k,\delta(\nabla_y\omega^k))\widetilde{H}^{\frac{1}{2}}\|^2_{L^2([0,u]\times\mathscr{O})} \leq$$
$$\frac{1}{2^k} + C_5\int_0^u\int_{x_0}^v H(y,x)\|(\delta\omega^{k-1},\nabla\delta\omega^{k-1})\|^2_{L^2(\mathscr{O})}dxdy\ .$$

This leads to the following inequalities:

$$\forall u\in[0,u_*],\ e^{-\Lambda u}\||x|^{-\alpha}(\delta\omega^k,\delta(\nabla_x\omega^k))(u)\|^2_{L^2([x_0,v[\times\mathscr{O})}$$
$$\leq \frac{1}{2^k} + C_5\int_0^u\int_{x_0}^v H(y,x)\|(\delta\omega^{k-1},\nabla\delta\omega^{k-1})\|^2_{L^2(\mathscr{O})}dxdy\ ; \tag{D.5}$$

$$\forall v\in[x_0,0[,\ |v|^{-\alpha}e^{-\Lambda v}\|(\delta\omega^k,\delta(\nabla_y\omega^k))(v)\|^2_{L^2([0,u]\times\mathscr{O})}$$
$$\leq \frac{1}{2^k} + C_5\int_0^u\int_{x_0}^v H(y,x)\|(\delta\omega^{k-1},\nabla\delta\omega^{k-1})\|^2_{L^2(\mathscr{O})}dxdy\ ; \tag{D.6}$$

$$\int_0^u\int_{x_0}^v H(y,x)\|(\delta\omega^k,\nabla\delta\omega^k)\|^2_{L^2(\mathscr{O})} \leq \frac{1}{2^k a(\Lambda)} + \sigma^2\int_0^u\int_{x_0}^v H(y,x)\|(\delta\omega^{k-1},\nabla\delta\omega^{k-1})\|^2_{L^2(\mathscr{O})}\ . \tag{D.7}$$

where

$$a(\Lambda) = 2\Lambda - C_5 - c_1(c_0,\bar{c}_0,n,h)\quad\text{and}\quad 0 < \sigma^2 < \frac{C_5}{a(\Lambda)} < \frac{1}{2}\ ,$$

provided that $\Lambda$ is large enough. For the last inequality, see D.2. $\square$

# E    Proof of corollary 3.16 page 16

Let $\varepsilon \in ]0, -x_0[$. Recall $\mathcal{D}_{*,\varepsilon} = [0, u_*] \times [x_0, -\varepsilon] \times \mathcal{O}$. In order to show that $\omega \in C^1(\mathcal{D}_*)$ we will show that $\omega \in C^1(\mathcal{D}_{*,\varepsilon})$ for any epsilon. We repeat what we did before to obtain that $\omega$ is continuous. Since

$$\sup_{k \in \mathbb{N}} \|\nabla \omega^k\|_{L^\infty(\mathcal{D}_{*,\varepsilon})} < C(C_0, x_0, \varepsilon) \;,$$

we only need to show that the sequence of second order derivatives $(\nabla^2 \omega^k)_{k \in \mathbb{N}}$ is bounded on $\mathcal{D}_{*,\varepsilon}$. This follows from (B.41), (B.42), (3.27) and (3.28) (recall $m - 2 > \frac{n-1}{2} + 2$). Thus again by Arzela-Ascoli theorem, the weak compactness and the interpolation theorem, one obtains that:

- the sequence (or a subsequence of it) $(\nabla \omega^k)_{k \in \mathbb{N}}$ converges uniformly towards $\nabla \omega$ on $\mathcal{D}_{*,\varepsilon}$,

- $\nabla \omega$ is a continuous function on $\mathcal{D}_{*,\varepsilon}$ and then on $\mathcal{D}_*$,

- $\forall s \in [0, m-2] \cap \mathbb{N}$, $\nabla \omega^k(y, x) \longrightarrow \nabla \omega(y, x)$ in $H^s(\mathcal{O})$ uniformly in $(y, x)$ on the compact $[0, u_*] \times [x_0, -\varepsilon]$ and that
$$\forall (y, x) \in [0, u_*] \times [x_0, 0[, \; \nabla \omega(u, v) \in C^2(\mathcal{O}) \;.$$

Let us show that $\omega \in C^2(\mathcal{D}_*)$. Again, we repeat the previous argument. Let $\varepsilon \in ]0, -x_0]$ be fixed. We already know that the sequence of second order derivatives $(\nabla^2 \omega^k)_{k \in \mathbb{N}}$ is uniformly bounded on $\mathcal{D}_{*,\varepsilon}$. Thus, it remains to show that the sequence of third order derivatives $(\nabla^3 \omega^k)_{k \in \mathbb{N}}$ is also uniformly bounded on $\mathcal{D}_{*,\varepsilon}$. From this property, it will follow that the sequence of second order derivatives is uniformly equicontinuous and then the theorem of Arzela-Ascoli applies. From some inequalities obtained so far, we see that the sequences $(\partial^3_{\mu\nu\beta} \omega^k)_{k \in \mathbb{N}}$ are uniformly bounded on $\mathcal{D}_{*,\varepsilon}$ for $\mu\nu\beta \neq xxx$ and $\mu\nu\beta \neq yyy$. Indeed we have:

- By choosing $j_0 \geq 3$ (which is the case since $m - 1 > 3 + \frac{n-1}{2}$), inequality (B.33) page 31 shows that the sequence $(\partial^3_{\mu\nu\beta} \omega^k)_{k \in \mathbb{N}}$ is uniformly bounded on $\mathcal{D}_{*,\varepsilon}$ for $\mu\nu\beta = xyA$.

- From inequalities (B.41) and (B.42) with $j_0 = 3$, we obtain that these sequences are uniformly bounded on $\mathcal{D}_{*,\varepsilon}$ for $\mu\nu\beta \in \{yyA, yAB, xxA, xAB\}$.

- Since $m > \frac{n+7}{2}$, the case $\mu\nu\beta = ABC$ will follows from inequality (3.27).

- The analysis of the right hand side of identity (B.36) gives the desired control in the case $\mu\nu\beta = xyy$ whereas $x\partial_x$-differentiating the partial differential equation satisfied by $\omega^{k+1}$ gives the result in the case $\mu\nu\beta = yxx$.

It thus remains to show that the sequences $(\partial^3_{yyy} \omega^k)_{k \in \mathbb{N}}$ and $(\partial^3_{xxx} \omega^k)_{k \in \mathbb{N}}$ are uniformly bounded on $\mathcal{D}_{*,\varepsilon}$. We start with $(\partial^3_{yyy} \omega^k)_{k \in \mathbb{N}}$. If we $\partial^2_y$-differentiate the differential equation satisfied by $\omega^{k+1}$ we obtain:

$$4\partial_x(\partial^3_y \omega^{k+1}) + \frac{n-1}{\rho} \partial^3_y \omega^{k+1} + (-x)^{-2} \partial^3_y \omega^k \left(\frac{\partial G}{\partial p_3}\right)^k (\ldots) = \Phi^k \tag{E.1}$$

where

$$\begin{aligned}\Phi^k = &= \frac{n-1}{\rho^2}(\partial_x - \partial_y)\partial_y\omega^{k+1}(y,x) + \frac{n-1}{2\rho^3}(\partial_x - \partial_y)\omega^{k+1} + \frac{n-1}{\rho}\partial_x\partial_y^2\omega^{k+1} \\ &+ \frac{h^{AB}}{\rho^2}\partial_A\partial_B\partial_y^2\omega^{k+1} + \frac{3h^{AB}}{2\rho^4}\partial_A\partial_B\omega^{k+1}(y,x) + \frac{2h^{AB}}{\rho^3}\partial_y\partial_A\partial_B\omega^{k+1}(y,x) \\ &- \frac{2\Gamma^B}{\rho^3}\partial_y\partial_B\omega^{k+1} - \frac{3\Gamma^B}{\rho^4}\partial_B\omega^{k+1} - \frac{\Gamma^B}{\rho^2}\partial_B\partial_y^2\omega^{k+1} \\ &- (-x)^{\frac{n+3}{2}}(\partial_y^2 G)^k(\ldots) - x^{-2}\partial_y^2\omega^k\left(\frac{\partial G}{\partial p_1}\right)^k(\ldots) - x^{-2}\partial_y\omega^k\partial_y\left(\left(\frac{\partial G}{\partial p_1}\right)^k(\ldots)\right) \\ &- x^{-2}\partial_y^2(x\partial_x\omega^k)\left(\frac{\partial G}{\partial p_2}\right)^k(\ldots) - x^{-2}\partial_y(x\partial_x\omega^k)\partial_y\left(\left(\frac{\partial G}{\partial p_2}\right)^k(\ldots)\right) \\ &- x^{-2}\partial_y^2\omega^k\partial_y\left(\left(\frac{\partial G}{\partial p_3}\right)^k(\ldots)\right) \\ &- x^{-2}\partial_y^2\partial_A\omega^k\left(\frac{\partial G}{\partial p_4}\right)^k(\ldots) + x^{-2}\partial_y\partial_A\omega^k\partial_y\left(\left(\frac{\partial G}{\partial p_4}\right)^k(\ldots)\right).\end{aligned}$$

From what has been said so far, we deduce that the coefficients of Equation (E.1) are uniformly bounded on $\mathcal{D}_{*,\varepsilon}$. Namely, we have $\|\frac{n-1}{\rho}\|_{L^\infty(\mathcal{D}_{*,\varepsilon})} < C$ and

$$\sup_{k\in\mathbb{N}}\|(-x)^{-2}\left(\frac{\partial G}{\partial p_3}\right)^k\|_{L^\infty(\mathcal{D}_{*,\varepsilon})} < C(x_0,\varepsilon),\ \sup_{k\in\mathbb{N}}\|\Phi^k\|_{L^\infty(\mathcal{D}_{*,\varepsilon})} < C(x_0,\varepsilon). \tag{E.2}$$

As we did before, we have:

$$\begin{aligned}\partial_x\left(\overline{H}(x)|\partial_y^3\omega^{k+1}|^2\right) &= (2\alpha|x|^{-1} - \Lambda)\overline{H}(x)|\partial_y^3\omega^{k+1}|^2 + 2\overline{H}(x)\partial_x\partial_y^3\omega^{k+1}.\partial_y^3\omega^{k+1} \\ &\leq -\Lambda\overline{H}(x)|\partial_y^3\omega^{k+1}|^2 + 2\overline{H}(x)\partial_x\partial_y^3\omega^{k+1}.\partial_y^3\omega^{k+1}.\end{aligned}$$

From (E.2), we deduce that

$$\begin{aligned}\partial_x\left(\overline{H}(x)|\partial_y^3\omega^{k+1}|^2\right) &\leq -\Lambda\overline{H}(x)|\partial_y\nabla_y\omega^{k+1}|^2 \\ &+ 2\overline{H}(x)\partial_y^3\omega^{k+1}\left(-\frac{n-1}{\rho}\partial_y^3\omega^{k+1} - (-x)^{-2}\partial_y^3\omega^k\left(\frac{\partial G}{\partial p_3}\right)^k(\ldots) + \Phi^k\right) \\ &\leq \overline{H}(x)\left((c + \frac{C}{\delta} - \Lambda)|\partial_y^3\omega^{k+1}|^2 + \delta C|\partial_y^3\omega^k|^2\right) + C|x|^{-\alpha}e^{-\Lambda x}.\end{aligned}$$

By choosing $\Lambda$ large enough, we have:

$$\partial_x\left(\overline{H}(x)|\partial_y^3\omega^{k+1}|^2\right) \leq \delta\overline{H}(x)C|\partial_y^3\omega^k|^2 + C|x|^{-\alpha}e^{-\Lambda x}, \tag{E.3}$$

which is then integrated in $x$ to obtain:

$$\overline{H}(x)|\partial_y^3\omega^{k+1}|^2(y,x) \leq \overline{H}(x_0)|\partial_y^3\omega^{k+1}|^2(y,x_0) + C\int_{x_0}^x e^{-\Lambda s}\left(\delta|s|^{-2\alpha}|\partial_y^3\omega^k|^2 + |s|^{-\alpha}\right)ds.$$

Equivalently this reads

$$\overline{H}(x)|\partial_y^3\omega^{k+1}|^2(y,x) \leq \overline{H}(x_0)|\partial_y^3\omega_0^{-,k+1}|^2(y,x_0) + C\int_{x_0}^x e^{-\Lambda s}\left(\delta|s|^{-2\alpha}|\partial_y^3\omega^k|^2 + |s|^{-\alpha}\right)ds$$

As we did many times before, for a convenient choice of δ, we obtain the estimate:

$$\begin{aligned}
|x|^{-2\alpha}|\partial_y^3\omega^{k+1}|^2(y,x) &\leq 2\left(C\int_{x_0}^0 |s|^{-\alpha}e^{-\Lambda s}ds + \sup_{k\in\mathbb{N},\ y\in[0,y_0]}\overline{H}(x_0)|\partial_y^3\omega_0^{-,k}|^2(y)\right)\\
&\leq c(x_0,\alpha)\left(1+\sup_{k\in\mathbb{N}}\|\omega_0^{-,k+1}\|^2_{C^3([0,y_0]\times\mathscr{O})}\right)\\
&\leq c(x_0,\alpha)\left(1+\sup_{k\in\mathbb{N}}\|\omega_0^{-,k+1}\|^2_{H^{m+1}(\mathcal{C}^-)}\right) < \infty
\end{aligned}$$

i.e.
$$|x|^{-\alpha}|\partial_y^3\omega^{k+1}|(y,x) < C\ . \tag{E.4}$$

The same holds for $|x|^{-\alpha}|\partial_x^3\omega^{k+1}|$ if instead we $\partial_x$–differentiate the differential equation satisfied by $\omega^{k+1}$: $|x|^{-\alpha}|\partial_x^3\omega^{k+1}|(y,x) < C$. Therefore we have proved that

$$\sup_{k\in\mathbb{N}}\|\nabla_y^3\omega^k\|_{L^\infty(\mathcal{D}_{*,\varepsilon})} < C\ . \tag{E.5}$$

It then follows that the family of functions $\{\nabla^2\omega^{k_j},\ j\in\mathbb{N}\}$ is uniformly equicontinuous and from the Arzela-Ascoli theorem there exists a subsequence of $(\nabla^2\omega^{k_j})_{j\in\mathbb{N}}$ denoted again by the same symbol which converges uniformly to a continuous function $\tilde\omega$ on $\mathcal{D}_{*,\varepsilon}$. Since the sequence $(\nabla\omega^{k_j})_{j\in\mathbb{N}}$ converges uniformly to $\nabla\omega$ on $\mathcal{D}_{*,\varepsilon}$ we conclude that $\nabla\omega$ is differentiable on $\mathcal{D}_{*,\varepsilon}$ and $\nabla^2\omega = \tilde\omega$ which proves that $\omega\in C^2(\mathcal{D}_{*,\varepsilon})$ for all $\varepsilon\in]0,-x_0]$ i.e. $\omega\in C^2(\mathcal{D}_*)$. To end the proof, it remains to show that $\omega$ solves the characteristic Cauchy problem (3.8). Recall that the partial differential equation satisfied by $\omega^{k+1}$ reads

$$\partial_x\partial_y\omega^{k+1} = \Psi^k \tag{E.6}$$

where

$$\begin{aligned}
\Psi^k &= \frac{(n-1)}{4\rho}\Big(\partial_x - \partial_y\Big)\omega^{k+1} + \frac{h^{AB}\partial_A\partial_B\omega^{k+1}}{4\rho^2}\\
&\quad -\frac{1}{4}\Gamma^B\partial_B\omega^{k+1} - \frac{1}{4}|x|^{-\frac{n+3}{2}}G(z,|x|^{\frac{n-1}{2}}(\omega^k,\nabla\omega^k))
\end{aligned}$$

being a continuous function on $\mathcal{D}_*$. To conclude we consider the limits pointwise in (E.6) and we are led to

$$\Box_{y,\eta}\omega = |x|^{-\frac{n+3}{2}}G(z,|x|^{\frac{n-1}{2}}(\omega,\nabla\omega))$$

thus $\omega$ is a classical solution of the characteristic initial value problem (3.8). □

**ACKNOWLEDGEMENTS:** RTW acknowledged useful discussions with Piotr Chruściel and is grateful to the group Gravitational Physics of Vienna University for hospitality during part of this work.

Figure 1: Characteristic cone $\mathcal{C}_{a,x}^+$ and its interior.

Figure 2: Images of the unbounded domain $\mathcal{Y}_{a,x}^+$ and the cone $\mathcal{C}_{a,x}^+$.

Figure 3: Neighborhood $V_{0,y}$ of the tip of the cone $\mathcal{C}_{-\frac{1}{a},y}^+$ and the cone $\mathcal{C}_{\lambda,y}^-$.

PSfrag replacements
$-\frac{1}{a}$
$\mathcal{C}^+$
$\mathcal{C}^-$
$V_{0,y}^\lambda$
$\lambda$

Figure 4: Truncated cones $\mathcal{C}^+$ et $\mathcal{C}^-$ .

PSfrag replacements
$\frac{1}{a}$
$-\frac{1}{a}$
$\mathcal{C}^+_{-\varepsilon,u}$
$\mathcal{C}^+$
$x$
$y$
$r$
$\rho$
$\lambda$
$\mathcal{C}^-$
$\mathcal{C}^-_{u,v}$
$\mathcal{C}^+_{u,v}$
$\mathcal{D}$

Figure 5: Future neighborhood $\mathcal{D}$ of the truncated cones $\mathcal{C}^+$ and $\mathcal{C}^-$ .